\input ./style/arxiv-general.cfg
\documentclass[aop,MSNbibl,seceqn,dvips]{arximspdf}
\makeatletter
   \@ifpackageloaded{graphicx}{}{\usepackage{graphicx}}
\makeatother
\usepackage{eufrak}
\usepackage{mathrsfs}

\doi{10.1214/14-AOP965} 
\volume{44}
\issue{1}
\pubyear{2016}
\firstpage{1}
\lastpage{41}
\docsubty{FLA}

\makeatletter
\def\geqslant{\geq}
\def\leqslant{\leq}
\newcommand{\rrvert}{\vert}
\newcommand{\rrVert}{\Vert}
\newcommand{\llvert}{\vert}
\newcommand{\llVert}{\Vert}
\newcommand{\eqref}[1]{(\ref{#1})}
\newtheorem{prop}{Proposition}[section]
\newtheorem{cor}[prop]{Corollary}
\newtheorem{lemme}[prop]{Lemma}
\newproclaim{rem}[prop]{Remark}
\newtheorem{teo}[prop]{Theorem}
\newproclaim{defi}[prop]{Definition}
\makeatother

\begin{document}
\begin{frontmatter}

\title{Quantitative stable limit theorems on the~Wiener~space}
\runtitle{Quantitative stable limit theorems on the Wiener space}

\begin{aug}
\author[A]{\fnms{Ivan}~\snm{Nourdin}\corref{}\ead[label=e1]{ivan.nourdin@uni.lu}\thanksref{T1}},
\author[B]{\fnms{David}~\snm{Nualart}\ead[label=e2]{nualart@math.ku.edu}\thanksref{T2}}
\and
\author[A]{\fnms{Giovanni}~\snm{Peccati}\ead[label=e3]{giovanni.peccati@gmail.com}\thanksref{T3}}
\runauthor{I. Nourdin, D. Nualart and G. Peccati}
\affiliation{Universit\'e de Lorraine, Kansas University and
Universit\'e du Luxembourg}
\address[A]{I. Nourdin\\
G. Peccati\\
FSTC---UR en Math\'{e}matiques\\
Universit\'e du Luxembourg\\
6 rue Richard Coudenhove-Kalergi\\
Luxembourg City 1359\\
Luxembourg\\
\printead{e1}\\
\phantom{E-mail:\ }\printead*{e3}} 
\address[B]{D. Nualart\\
Department of Mathematics\\
University of Kansas\\
1460 Jayhawk Boulevard\\
Lawrence, Kansas 66045\\
USA\\
\printead{e2}}
\end{aug}
\thankstext{T1}{Supported in part by the
French ANR Grant ANR-10-BLAN-0121.}
\thankstext{T2}{Supported in part by NSF Grant
DMS-12-08625.}
\thankstext{T3}{Supported in part by Grant F1R-MTH-PUL-12PAMP
(PAMPAS), from
Luxembourg University.}

\received{\smonth{5} \syear{2013}}
\revised{\smonth{8} \syear{2014}}

%
\begin{abstract}
We use Malliavin operators in order
to prove quantitative stable limit theorems on the Wiener space, where
the target distribution is given by a possibly multidimensional mixture
of Gaussian distributions. Our findings refine and generalize previous
works by Nourdin and Nualart [\textit{J. Theoret. Probab.} \textbf{23} (2010) 39--64]
and Harnett and Nualart
[\textit{Stochastic Process. Appl.} \textbf{122} (2012) 3460--3505], and
provide a substantial contribution to a recent line of research,
focussing on limit theorems on the Wiener space, obtained by means of
the Malliavin calculus of variations. Applications are given to
quadratic functionals and weighted quadratic variations of a fractional
Brownian motion.
\end{abstract}

%
\end{frontmatter}

\section{Introduction and overview}\label{sec1}

Originally introduced by R\'{e}nyi in the landmark paper \cite{reny},
the notion of \textit{stable convergence} for random variables (see
Definition~\ref{dstable} below) is an intermediate concept, bridging
convergence in distribution (which is a weaker notion) and convergence
in probability (which is stronger). One crucial feature of stably
converging sequences is that they can be naturally paired with
sequences converging in probability (see, e.g., the statement of Lemma~\ref{lstable} below), thus yielding a vast array of noncentral limit
results---most notably convergence toward \textit{mixtures} of Gaussian
distributions. This last feature makes indeed stable convergence
extremely useful for applications, in particular to the asymptotic
analysis of functionals of semimartingales, such as power variations,
empirical covariances, and other objects of statistical relevance. See
the classical reference \cite{JacSh}, Chapter VIII.5, as well as the
recent survey \cite{podvetter}, for a discussion of stable convergence
results in a semimartingale context.

Outside the (semi)martingale setting, the problem of characterizing
stably converging sequences is for the time being much more delicate.
Within the framework of limit theorems for functionals of general
Gaussian fields, a step in this direction appears in the paper \cite{pectud}, by Peccati and Tudor, where it is shown that central limit
theorems (CLTs) involving sequences of multiple Wiener--It\^o integrals
of order${}\geqslant{}2$ are always stable. Such a result is indeed an
immediate consequence of a general multidimensional CLT for chaotic
random variables, and of the well-known fact that the first Wiener
chaos of a Gaussian field coincides with the $L^2$-closed Gaussian
space generated by the field itself (see \cite{np-book}, Chapter~6, for
a general discussion of multidimensional CLTs on the Wiener space).
Some distinguished applications of the results in \cite{pectud} appear,
for example, in the two papers \cite{cnw,BNCP}, respectively, by
Corcuera et al.  and by Barndorff-Nielsen et al., where the
authors establish stable limit theorems (toward a Gaussian mixture) for
the power variations of pathwise stochastic integrals with respect to a
Gaussian process with stationary increments. See \cite{noupecwei} for
applications to the weighted variations of an iterated Brownian motion.
See \cite{boupec} for some quantitative analogues of the findings of
\cite{pectud} for functionals of a Poisson measure.

Albeit useful for many applications, the results proved in \cite{pectud} do not provide any intrinsic criterion for stable convergence
toward Gaussian mixtures. In particular, the applications developed in
\cite{BNCP,cnw,noupecwei} basically require that one is able to
represent a given sequence of functionals as the combination of three
components---one converging in probability to some nontrivial random
element, one living in a~finite sum of Wiener chaoses and one vanishing
in the limit---so that the results from \cite{pectud} can be directly
applied. This is in general a highly nontrivial task, and such a
strategy is technically too demanding to be put into practice in
several situations (e.g., when the chaotic decomposition of a
given functional cannot be easily computed or assessed).

The problem of finding effective intrinsic criteria for stable
convergence on the Wiener space toward mixtures of Gaussian
distributions---without resorting to chaotic decompositions---was
eventually tackled by Nourdin and Nualart in \cite{nounua}, where one
can find general sufficient conditions ensuring that a sequence of \textit{multiple Skorohod integrals} stably converges to a mixture of Gaussian
distributions. Multiple Skorohod integrals are a generalization of
multiple Wiener--It\^o integrals (in particular, they allow for random
integrands), and are formally defined in Section~\ref{ssmall} below.
It is interesting to note that the main results of \cite{nounua} are
proved by using a generalization of a characteristic function method,
originally applied by Nualart and Ortiz-Latorre in \cite{nuaort} to
provide a Malliavin calculus proof of the CLTs established in \cite{nunugio,pectud}.
In particular, when specialized to multiple
Wiener--It\^o integrals, the results of \cite{nounua} allow to recover
the ``fourth moment theorem'' by Nualart and Peccati \cite{nunugio}. A
first application of these stable limit theorems appears in \cite{nounua}, Section~5, where one can find stable mixed Gaussian limit
theorems for the weighted quadratic variations of the fractional
Brownian motion (fBm), complementing some previous findings from \cite{nounuatud}.
Another class of remarkable applications of the results of
\cite{nounua} are the so-called \textit{It\^o formulae in law}; see \cite{HarNua1,HarNua2,nourev,nourevswa}. Reference \cite{HarNua1} also
contains some multidimensional extensions of the abstract results
proved in \cite{nounua} (with a proof again based on the characteristic
function method). Further applications of these techniques can be found
in \cite{rev2009}. An alternative approach to stable convergence on the
Wiener space, based on decoupling techniques, has been developed by
Peccati and Taqqu in \cite{petamwii}.

One evident limitation of the abstract results of \cite{HarNua1,nounua} is that they do not provide any information about rates of
convergence. The aim of this paper is to prove several \textit{quantitative versions} of the abstract results proved in
\cite{HarNua1,nounua}, that is, statements allowing one to explicitly assess
quantities of the type
\[
\bigl| E \bigl[\varphi\bigl(\delta^{q_1}(u_1),\ldots,
\delta^{q_d}(u_d)\bigr)\bigr] - E\bigl[\varphi(F) \bigr] \bigr|,
\]
where $\varphi$ is an appropriate test function on $\mathbb{R}^d$,
each $\delta^{q_i}(u_i)$ is a multiple Skorohod integral of order $q_i\geqslant1$, and
$F$ is a $d$-dimensional mixture of Gaussian distributions. Most
importantly, we shall show that our bounds also yield natural
sufficient conditions for stable convergence toward $F$. To do this, we
must overcome a number of technical difficulties, in particular:
\begin{itemize}
\item We will work in a general framework and without any
underlying semimartingale structure, in such a way that the powerful
theory of stable convergence for semimartingales (see again \cite{JacSh}) cannot be applied.

\item Although there are many versions of Stein's method allowing
one to deal with general continuous non-Gaussian targets (see, e.g.,
\cite{chatterjee,D,EdVie,EdViq,KT1,KT2,R}), it seems that none of
them can be reasonably applied to the limit theorems that are studied
in this paper. Indeed, the above quoted contributions fall mainly in
two categories: either those requiring that the density of the target
distribution is explicitly known (and in this case the so-called
``density approach'' can be applied---see, e.g., \cite{chatterjee,D,EdVie,EdViq}), or those requiring that the target distribution is the
invariant measure of some diffusion process (so that the ``generator
approach'' can be used---see, e.g., \cite{KT1,KT2,R}). In both
instances, a detailed analytical description of the target distribution
must be available. In contrast, in the present paper we consider limit
distributions given by the law of random elements of the type $S\cdot
\eta=(S_1 \eta_1,\ldots,S_d\eta_d)$, where $\eta=(\eta_1,\ldots,\eta
_d)$ is
a Gaussian vector, and $S =(S_1,\ldots,S_d)$ is an independent random
element that is suitably regular in the sense of Malliavin calculus. In
particular, in our framework \textit{no a priori knowledge} of the
distribution of $S$ (and therefore of $S\cdot\eta$) is required. One
should note that in \cite{chatterjee} one can find an application of
Stein's method to the law of random objects with the form $S\eta$,
where $\eta$ is a one-dimensional Gaussian random variable and $S$ has
a law with a two-point support (of course, in this case the density of
$S\eta$ can be directly computed by elementary arguments).
\end{itemize}

Our techniques rely on an interpolation procedure and on the use of
Malliavin operators. To our knowledge, the main bounds proved in this
paper, that is, the ones appearing in Proposition~\ref{prop1-1},
Theorems~\ref{tkol1} and~\ref{tmain}, are first ever explicit
upper bounds for mixed normal approximations in a nonsemimartingale setting.

Note that, in our discussion, we shall separate the case of
one-dimensional Skorohod integrals of order 1 (discussed in Section~\ref{sdq1}) from the general case (discussed in Section~\ref{sdqany}),
since in the former setting one can exploit some useful
simplifications, as well as obtain some effective bounds in the
Wasserstein and Kolmogorov distances. As discussed below, our results
can be seen as abstract versions of classic limit theorems for Brownian
martingales, such as the ones discussed in \cite{RY}, Chapter~VIII.

Although our results deal only with Skorohod integrals, they can be
applied in the context of Stratonovich integrals. In fact, the
Stratonovich integral can be expressed as a Skorohod integral plus a
complementary term and in many problems this complementary term does
not contribute to the limit. Examples of this situation are the It\^o
formulas in law for different types of Stratonovich integrals obtained
by Harnett and Nualart in \cite{HarNua1,HarNua2} and the weak
convergence of weighted variations established by Nourdin and Nualart
in \cite{nounua}.

To illustrate our findings, we provide applications to quadratic
functionals of a fractional Brownian motion (Section~\ref{ssex1}) and
to weighted quadratic variations (Section~\ref{sex2}). The results of
Section~\ref{ssex1} generalize some previous findings by Peccati and
Yor \cite{py1,py2}, whereas those of Section~\ref{sex2} complement
some findings by Nourdin, Nualart and Tudor \cite{nounuatud}.

The paper is organized as follows. Section~\ref{sec2} contains some
preliminaries on Gaussian analysis and stable convergence.
In Section~\ref{sdq1}, we first derive estimates for the distance between the laws of a
Skorohod integral of order $1$ and of a mixture of Gaussian
distributions (see Proposition~\ref{prop1-1}). As a corollary, we
deduce the stable limit theorem for a sequence of multiple Skorohod
integrals of order $1$ obtained in \cite{HarNua1}, and we obtain rates
of convergence in the Wasserstein and Kolmogorov distances. We apply
these results to a sequence of quadratic functionals of the fractional
Brownian motion. Section~\ref{sec4} contains some additional notation and a
technical lemma that are used in Section~\ref{sdqany} to establish bounds in the
multidimensional case for Skorohod integrals of general orders.
Finally, in Section~\ref{sex2} we present the applications of these results to
the case of weighted quadratic variations of the fractional Brownian
motion. The \hyperref[sec7]{Appendix} contains some technical lemmas needed in Section~\ref{sex2}.

\section{Gaussian analysis and stable convergence}\label{sec2}

In the next two subsections, we discuss some basic notions of Gaussian
analysis and Malliavin calculus. The reader is referred to the
monographs \cite{nualartbook} and \cite{np-book} for any unexplained
definition or result.

\subsection{Elements of Gaussian analysis}\label{ssmall}

Let $\EuFrak{H}$ be a real separable infinite-dimensional Hilbert space.
For any integer $q\geqslant1$, we denote by $\EuFrak{H}^{\otimes q}$ and $\EuFrak{H}^{\odot q}$, respectively, the
$q$th tensor product and the $q$th symmetric tensor product of $\EuFrak{H}$.
In what follows, we write $X=\{X(h)\dvtx  h\in\EuFrak{H}\}$ to
indicate an isonormal Gaussian process over
$\EuFrak{H}$. This means that $X$ is a centered Gaussian family, defined
on some probability space $(\Omega,\mathcal{F},P)$, with a covariance
structure given by
%
\begin{equation}
\label{v1}
E\bigl[X(h)X(g)\bigr] = \langle h,g \rangle_\EuFrak{H},
\qquad h,g\in\EuFrak{H}.
\end{equation}
From now on, we assume that $\mathcal{F}$ is the $P$-completion of the
$\sigma$-field generated by $X$. For every integer $q\geqslant1$, we let
$\mathcal{H}_{q}$ be the $q$th \textit{Wiener chaos} of $X$,
that is, the closed linear subspace of $L^{2}(\Omega)$
generated by the random variables $\{H_{q}(X(h)),h\in\EuFrak{H},\llVert
h\rrVert_{\EuFrak{H}}=1\}$, where $H_{q}$ is the $q$th Hermite
polynomial defined by
\[
H_q(x)=(-1)^q e^{x^2/2}\,\frac{d^q}{dx^q}
\bigl(e^{-x^2/2} \bigr).
\]
We denote by $\mathcal{H}_{0}$ the space of constant random variables. For
any $q\geqslant1$, the mapping $I_{q}(h^{\otimes q})=q!H_{q}(X(h))$
provides a
linear isometry between $\EuFrak{H}^{\odot q}$
(equipped with the modified norm $\sqrt{q!}\llVert  \cdot\rrVert
_{\EuFrak{H}^{\otimes q}}$) and $\mathcal{H}_{q}$ [equipped with the $L^2(\Omega)$
norm]. For $q=0$, we set by convention $\mathcal{H}_{0}=\mathbb{R}$ and $I_{0}$ equal to the identity map.

It is well known (Wiener chaos expansion) that $L^{2}(\Omega)$
can be decomposed into the infinite orthogonal sum of the spaces
$\mathcal{H}_{q}$, that is: any square integrable random variable $F\in
L^{2}(\Omega)$ admits the following chaotic expansion:
%
\begin{equation}\label{E}
F=\sum_{q=0}^{\infty}I_{q}(f_{q}),
\end{equation}
where $f_{0}=E[F]$, and the $f_{q}\in\EuFrak{H}^{\odot q}$, $q\geqslant
1$, are
uniquely determined by $F$. For every $q\geqslant0$, we denote\vspace*{1pt} by
$J_{q}$ the
orthogonal projection operator on the $q$th Wiener chaos. In
particular, if $F\in L^{2}(\Omega)$ is as in (\ref{E}), then $J_{q}F=I_{q}(f_{q})$ for every $q\geqslant0$.

Let $\{e_{k}, k\geqslant1\}$ be a complete orthonormal system in
$\EuFrak H$.
Given $f\in\EuFrak{H}^{\odot p}$, $g\in\EuFrak H^{\odot q}$ and
$r\in\{0,\ldots,p\wedge q\}$, the $r$th \textit{contraction} of $f$ and $g$
is the element of $\EuFrak{H}^{\otimes(p+q-2r)}$ defined by
%
\begin{equation} \label{v2}
f\otimes_{r}g=\sum_{i_{1},\ldots,i_{r}=1}^{\infty}
\langle f,e_{i_{1}}\otimes\cdots\otimes e_{i_{r}}
\rangle_{\EuFrak H^{\otimes
r}}\otimes\langle g,e_{i_{1}}\otimes\cdots\otimes
e_{i_{r}}\rangle_{\EuFrak{H}^{\otimes r}}.
\end{equation}
Notice that $f\otimes_{r}g$ is not necessarily symmetric. We denote its
symmetrization by $f\,\widetilde{\otimes}_{r}\,g\in\EuFrak{H}^{\odot
(p+q-2r)}$.
Moreover, $f\otimes_{0}g=f\otimes g$ equals the tensor product of
$f$ and
$g$ while, for $p=q$, $f\otimes_{q}g=\langle f,g\rangle_{\EuFrak{H}^{\otimes q}}$.
Contraction operators are useful for dealing with
products of multiple Wiener--It\^o integrals.

In the particular case where $\EuFrak{H}=L^{2}(A,\mathcal{A},\mu)$,
with $(A,\mathcal{A})$ is a measurable space and $\mu$ is a $\sigma$-finite and
nonatomic measure, one has that $\EuFrak{H}^{\odot q}=L_{s}^{2}(A^{q},
\mathcal{A}^{\otimes q},\mu^{\otimes q})$ is the space of symmetric and
square integrable functions on~$A^{q}$. Moreover, for every $f\in
\EuFrak{H}^{\odot q}$, $I_{q}(f)$ coincides with the multiple Wiener--It\^{o} integral
of order $q$ of $f$ with respect to $X$ (as defined, e.g., in \cite{nualartbook}, Section~1.1.2)
and (\ref{v2}) can be written as
\begin{eqnarray*}
&&(f\otimes_{r}g) (t_1,\ldots,t_{p+q-2r}) \\
&& \qquad =\int
_{A^{r}}f(t_{1},\ldots,t_{p-r},s_{1},
\ldots,s_{r})\\
&& \qquad \hspace*{26pt}{}\times
g(t_{p-r+1},\ldots,t_{p+q-2r},s_{1},
\ldots ,s_{r})\,d\mu (s_{1})\cdots \,d\mu(s_{r}).
\end{eqnarray*}

\subsection{Malliavin calculus}

Let us now introduce some elements of the Malliavin calculus of
variations with respect
to the isonormal Gaussian process $X$. Let $\mathcal{S}$
be the set of all smooth and cylindrical random variables of
the form
%
\begin{equation}\label{v3}
F=g \bigl( X(\phi_{1}),\ldots,X(\phi_{n}) \bigr),
\end{equation}
where $n\geqslant1$, $g\dvtx \mathbb{R}^{n}\rightarrow\mathbb{R}$ is a infinitely
differentiable function with compact support, and $\phi_{i}\in\EuFrak{H}$.
The \textit{Malliavin derivative} of $F$ with respect to $X$ is the
element of $L^{2}(\Omega,\EuFrak{H})$ defined as
\[
DF = \sum_{i=1}^{n}\frac{\partial g}{\partial x_{i}}
\bigl( X(\phi _{1}),\ldots,X(\phi_{n}) \bigr)
\phi_{i}.
\]
By iteration, one can
define the $q$th derivative $D^{q}F$ for every $q\geqslant2$, which is an
element of $L^{2}(\Omega,\EuFrak{H}^{\odot q})$.

For $q\geqslant1$ and $p\geqslant1$, ${\mathbb{D}}^{q,p}$ denotes
the closure
of $\mathcal{S}$ with respect to the norm $\Vert\cdot\Vert_{\mathbb
{D}^{q,p}}$, defined by
the relation
\[
\Vert F\Vert_{\mathbb{D}^{q,p}}^{p} = E \bigl[ |F|^{p} \bigr]
+\sum_{i=1}^{q}E \bigl( \bigl\Vert
D^{i}F\bigr\Vert_{\EuFrak{H}^{\otimes i}}^{p} \bigr).
\]
The Malliavin derivative $D$ verifies the following chain rule.
If $\varphi\dvtx \mathbb{R}^{n}\rightarrow\mathbb{R}$ is continuously
differentiable with bounded partial derivatives and if $F=(F_{1},\ldots
,F_{n})$ is a vector of elements of ${\mathbb{D}}^{1,2}$, then
$\varphi
(F)\in{\mathbb{D}}^{1,2}$ and
\[
D\varphi(F)=\sum_{i=1}^{n}
\frac{\partial\varphi}{\partial x_{i}}(F)DF_{i}.
\]

We denote by $\delta$ the adjoint of the operator $D$, also called the
\textit{divergence operator} or \textit{Skorohod integral} (see, e.g., \cite{nualartbook}, Section~1.3.2,  for an explanation of this terminology).
A random element $u\in L^{2}(\Omega,\EuFrak{H})$ belongs to the domain of $\delta$, noted
$\operatorname{Dom}\delta$,
if and
only if it verifies
\[
\bigl|E \bigl(\langle DF,u\rangle_{\EuFrak{H}} \bigr) \bigr|\leqslant c_{u}
\sqrt{E\bigl(F^2\bigr)}
\]
for any $F\in\mathbb{D}^{1,2}$, where $c_{u}$ is a constant depending only
on $u$. If $u\in\operatorname{Dom}\delta$, then the random variable
$\delta(u)$
is defined by the duality relationship (called ``integration by parts
formula''):
%
\begin{equation}\label{ipp}
E\bigl(F\delta(u)\bigr)=E \bigl(\langle DF,u\rangle_{\EuFrak H} \bigr),
\end{equation}
which holds for every $F\in{\mathbb{D}}^{1,2}$. The formula (\ref%
{ipp}) extends to the multiple Skorohod integral $\delta^{q}$, and we
have
%
\begin{equation}\label{dual}
E \bigl( F\delta^{q}(u) \bigr) =E \bigl( \bigl\langle
D^{q}F,u \bigr\rangle_{%
\EuFrak H^{\otimes q}} \bigr),
\end{equation}
for any element $u$ in the domain of $\delta^{q}$ and any random
variable $F\in\mathbb{D}^{q,2}$. Moreover, $\delta^{q}(h)=I_{q}(h)$ for any
$h\in\EuFrak{H}^{\odot q}$.

The following statement will be used in the paper, and is proved in
\cite{nounua}.

\begin{lemme}\label{lemme}
Let $q\geqslant1$ be an integer. Suppose that $F\in{\mathbb
{D}}^{q,2}$, and
let~$u$ be a symmetric element in $\operatorname{Dom}\delta^{q}$. Assume that,
for any $ 0\leqslant r+j\leqslant q$,
$ \langle D^{r}F,\delta^{j}(u) \rangle_{\EuFrak{H}^{\otimes
r}}\in
L^{2}(\Omega,\EuFrak{H}^{\otimes q-r-j})$. Then, for any $r=0,\ldots
,q-1$, $ \langle D^{r}F,u \rangle_{\EuFrak{H}^{\otimes r}}$ belongs
to the
domain of $\delta^{q-r}$ and we have
%
\begin{equation}\label{t3}
F\delta^{q}(u)=\sum_{r=0}^{q}
\pmatrix{q \cr
r}\delta^{q-r} \bigl( \bigl\langle D^{r}F,u
\bigr\rangle_{\EuFrak H^{\otimes r}} \bigr).
\end{equation}
[With the convention that $\delta^0(v)=v$, $v\in L^2(\Omega)$ and
$D^0F=F$, $F\in L^2(\Omega)$.]
\end{lemme}

For any Hilbert space $V$, we denote by $\mathbb{D}^{k,p}(V)$ the
corresponding Sobolev space of $V$-valued random variables (see \cite{nualartbook}, page~31).
The operator $\delta^q$ is continuous from $\mathbb{D}^{k,p}(\EuFrak
H^{\otimes q})$
to $\mathbb{D}^{k-q,p}$, for any $p>1$ and any integers $k\ge q\ge1$,
that is,
we have
%
\begin{equation}\label{Me2}
\bigl\llVert \delta^q (u)\bigr\rrVert _{\mathbb{D}^{k-q,p}}\leqslant
c_{k,p}\llVert u\rrVert _{\mathbb{D}^{k,p}(\EuFrak H^{\otimes q})},
\end{equation}
for all $u\in\mathbb{D}^{k,p}(\EuFrak H^{\otimes q})$, and some
constant $c_{k,p}>0$.
These estimates are consequences of Meyer inequalities (see \cite{nualartbook}, Proposition~1.5.7).
In particular, these estimates imply that
$\mathbb{D}^{q,2}(\EuFrak H^{\otimes q})\subset\operatorname{Dom}\delta
^{q}$ for any integer $q\geqslant1$.

The following commutation relationship between the
Malliavin derivative and the Skorohod integral (see \cite{nualartbook}, Proposition~1.3.2) is also useful:
%
\begin{equation}\label{comm1}
D\delta(u)=u+\delta(Du),
\end{equation}
for any $u\in\mathbb{D}^{2,2}(\EuFrak H)$. By induction, we can show
the following formula for any symmetric element $u$ in $\mathbb
{D}^{j+k,2}(\EuFrak{H}^{\otimes j})$
%
\begin{equation}\label{t5}
D^{k}\delta^{j}(u)=\sum_{i=0}^{j\wedge k}
\pmatrix{k
\cr
i}\pmatrix{j
\cr
i}i!\delta ^{j-i}\bigl(D^{k-i}u
\bigr).
\end{equation}
Also, we will make sometimes use of the following formula for the
variance of a~multiple
Skorohod integral. Let $u,v\in\mathbb{D}^{2q,2}(\EuFrak{H}^{\otimes
q})\subset \operatorname{Dom}\delta^q$
be two symmetric functions. Then
%
\begin{eqnarray}
E\bigl(\delta^{q}(u)\delta^{q}(v)\bigr) &=&E\bigl( \bigl
\langle u,D^{q}\bigl(\delta ^{q}(v)\bigr) \bigr
\rangle_{\EuFrak H^{\otimes q}}\bigr) \nonumber\\
\label{t13} &=&  \sum_{i=0}^{q}
\pmatrix{q
\cr
i}^{2}i!E \bigl( \bigl\langle u,\delta ^{q-i}
\bigl(D^{q-i}v\bigr) \bigr\rangle_{\EuFrak H^{\otimes q}} \bigr)
\\
\nonumber
&=&\sum_{i=0}^{q}\pmatrix{q
\cr
i}^{2}i!E \bigl( D^{q-i}u \,\widehat{\otimes
}_{2q-i}\,D^{q-i}v \bigr),
\end{eqnarray}
with the notation
\begin{eqnarray*}
&& D^{q-i}u\, \widehat{\otimes}_{2q-i}\,D^{q-i}v\\
&& \qquad  =\sum
_{ j,k,\ell=1} ^\infty \bigl\langle D^{q-i}
\langle u, \xi_j \otimes\eta_\ell\rangle_{\EuFrak{H}^{\otimes q}},
\xi_k \bigr\rangle _{\EuFrak H^{\otimes q-i}} \bigl\langle D^{q-i}
\langle v , \xi_k \otimes\eta_\ell\rangle_{\EuFrak{H}^{\otimes q}}, \xi_j \bigr\rangle_{\EuFrak H^{\otimes q-i}},
\end{eqnarray*}
where $\{\xi_j, j\geqslant1\}$ and $\{\eta_\ell, \ell\geqslant1\}$ are complete
orthonormal systems in $\EuFrak{H}^{\otimes q-i}$ and
$\EuFrak{H}^{\otimes i}$, respectively.

The operator $L$ is defined on the Wiener chaos expansion as
$L=\sum_{q=0}^{\infty}-qJ_{q}$,
and is called the \textit{infinitesimal generator of the Ornstein--Uhlenbeck
semigroup}. The domain of this operator in $L^{2}(\Omega)$ is the set
\[
\operatorname{Dom}L=\Biggl\{F\in L^{2}(\Omega)\dvtx \sum
_{q=1}^{\infty}q^{2}\llVert J_{q}F
\rrVert _{L^{2}(\Omega)}^{2}<\infty\Biggr\}=\mathbb{D}^{2,2}.
\]
There is an important relationship between the operators $D$, $\delta$
and $L$
(see \cite{nualartbook}, Proposition~1.4.3). A random variable $F$
belongs to the
domain of $L$ if and only if $F\in\operatorname{Dom} ( \delta D
) $
(i.e., $F\in{\mathbb{D}}^{1,2}$ and $DF\in\operatorname{Dom}\delta$), and in
this case
%
\begin{equation}\label{k1}
\delta DF=-LF.
\end{equation}
Note also that a random variable $F$ as in (\ref{E}) is in ${\mathbb{D}}
^{1,2}$ if and only if
$\sum_{q=1}^{\infty}qq!\Vert f_{q}\Vert_{\EuFrak H^{\otimes
q}}^{2}<\infty$,
and, in this case, $E ( \Vert DF\Vert_{\EuFrak H}^{2} )
=\sum_{q\geqslant1}qq!\Vert f_{q}\Vert_{\EuFrak H^{\otimes q}}^{2}$. If
$\EuFrak H=%
L^{2}(A,\mathcal{A},\mu)$ (with $\mu$ nonatomic), then the
derivative of a random variable $F$ as in (\ref{E}) can be identified with
the element of $L^{2}(A\times\Omega)$ given by
%
\begin{equation}\label{dtf}
D_{a}F=\sum_{q=1}^{\infty}qI_{q-1}
\bigl( f_{q}(\cdot,a) \bigr),\qquad a\in A.
\end{equation}

\subsection{Stable convergence}

The notion of stable convergence used in this paper is provided in the
next definition. Recall that the probability space $(\Omega, \mathcal
{F}, P)$ is such that $\mathcal{F}$ is the $P$-completion of the
$\sigma$-field generated by the isonormal process~$X$.

\begin{defi}[(Stable convergence)]\label{dstable}
Fix $d\geqslant1$. Let $\{F_{n}\}$ be a sequence of random variables with
values in $\mathbb{R}^{d}$, all defined on the probability space
$(\Omega, \mathcal{F}, P)$. Let $F$ be a $\mathbb{R}^d$-valued random variable
defined on
some extended probability space $(\Omega', \mathcal{F}', P')$. We say
that $F_{n}$ \textit{converges stably} to $F$, written $F_n \stackrel{\mathrm{st}}{\rightarrow}F$, if
%
\begin{equation}\label{estable}
\lim_{n \rightarrow\infty}E \bigl[Ze^{i \langle
\lambda, F_{n}  \rangle_{\mathbb{R}^{d}} } \bigr] = E'
\bigl[Ze^{i \langle
\lambda, F  \rangle_{\mathbb{R}^{d}} }\bigr],
\end{equation}
for every $\lambda\in\mathbb{R}^{d}$ and every bounded $\mathcal{F}$-measurable random variable $Z$.
\end{defi}

Choosing $Z=1$ in (\ref{estable}), we see that stable convergence
implies convergence in distribution. For future reference, we now list
some useful properties of stable convergence. The reader is referred,
for example, to \cite{JacSh}, Chapter~4,  for proofs. From now on, we
will use the symbol $\stackrel{P}{\rightarrow}$ to indicate
convergence in probability with respect to $P$.

\begin{lemme}\label{lstable}
Let $d\geqslant1$, and let $\{F_n\}$ be a
sequence of random variables with values in $\mathbb{R}^d$.
\begin{longlist}[4.]
\item[1.] $F_{n} \stackrel{\mathrm{st}}{\rightarrow}F$ if and only if $(F_n,
Z)\stackrel{\mathrm{law}}{\rightarrow}(F,Z)$, for every $\mathcal{F}$-measurable random
variable $Z$.

\item[2.] $F_{n} \stackrel{\mathrm{st}}{\rightarrow}F$ if and only if $(F_n,
Z)\stackrel{\mathrm{law}}{\rightarrow}(F,Z)$, for every random variable $Z$ belonging to
some set $\mathscr{Z} = \{Z_\alpha\dvtx  \alpha\in A\}$ such that the
$P$-completion of $\sigma(\mathscr{Z})$ coincides with $ \mathcal{F}$.

\item[3.] If $F_{n} \stackrel{\mathrm{st}}{\rightarrow}F$ and $F$ is $\mathcal{F}$-measurable, then
necessarily $F_n\stackrel{P}{\rightarrow}F$.

\item[4.] If $F_{n} \stackrel{\mathrm{st}}{\rightarrow}F$ and $\{Y_n\}$ is another
sequence of random
elements, defined on $(\Omega, \mathcal{F}, P)$ and such that $Y_n
\stackrel{P}{\rightarrow} Y$, then $(F_n, Y_n) \stackrel{\mathrm{st}}{\rightarrow}(F,Y)$.
\end{longlist}
\end{lemme}

The following statement (to which we will compare many results of the
present paper) contains criteria for the stable convergence of vectors
of multiple Skorohod integrals of the same order. The case $d=1$ was
proved in \cite{nounua}, Corollary~3.3, whereas the case of a general
$d$ is dealt with in \cite{HarNua1}, Theorem~3.2. Given $d\geqslant1$,
$\bolds{\mu} \in\mathbb{R}^d$ and a nonnegative definite $d\times d$
matrix $C$,
we shall denote by $\mathcal{N}_d(\bolds{\mu}, C)$ the law of a
$d$-dimensional Gaussian vector with mean $\bolds{\mu}$ and covariance
matrix~$C$.

\begin{teo} \label{tHN}
Let $q ,d \geqslant1$ be integers, and suppose that $F_n$ is a
sequence of
random variables in $\mathbb{R}^d$ of the form $F_n = \delta^q (u_n) =
 ( \delta^q (u_n^1), \ldots,\break  \delta^q(u_n^d)  )$, for a~sequence of $\mathbb{R}^d$-valued symmetric functions $u_n$ in
$\mathbb{D}^{2q, 2q}(\mathfrak{H}^{\otimes q})$. Suppose that the sequence
$F_n$ is bounded in $L^1(\Omega)$ and that:
\begin{enumerate}
\item $\langle u_n^j, \bigotimes_{\ell=1}^m (D^{a_\ell}F_n^{j_\ell})
\otimes h \rangle_{\mathfrak{H}^{\otimes q}}$ converges to zero in
$L^1(\Omega)$ for all integers $1 \leqslant j, j_\ell\leqslant d$,
all integers
$1 \leqslant a_1, \ldots, a_m, r \leqslant q-1$ such that $a_1 + \cdots
+ a_m + r
= q$, and all $h \in\mathfrak{H}^{\otimes r}$.
\item For each $1 \leqslant i,j \leqslant d$, $\langle  u_n^i,
D^qF_n^j\rangle
_{\mathfrak{H}^{\otimes q}}$ converges in $L^1(\Omega)$ to a random
variable $s_{ij}$, such that the random matrix $\Sigma:=  ( s_{ij}
 )_{d \times d}$ is nonnegative definite.
\end{enumerate}
Then $F_n \stackrel{\mathrm{st}}{\rightarrow}F$, where $F$ is a random variable
with values in
${\mathbb R}^d$ and with conditional Gaussian distribution ${\mathcal N}_d
(0, \Sigma)$ given $X$.
\end{teo}

\subsection{Distances}

For future reference, we recall the definition of some useful distances
between the laws of two real-valued random variables $F,G$.
\begin{itemize}
\item The \textit{Wasserstein distance} between the laws of $F$ and
$G$ is defined by
\[
d_{\mathrm{W}}(F,G)= \sup_{\varphi\in \operatorname{Lip} (1)} \bigl|E\bigl[ \varphi(F)\bigr]
-E\bigl[\varphi(G)\bigr]\bigr|,
\]
where $\operatorname{Lip}(1)$ indicates the collection of all Lipschitz
functions $\varphi$ with Lipschitz constant less than or equal to $1$.

\item The \textit{Kolmogorov distance} is
\[
d_{\mathrm{Kol}}(F,G)= \sup_{x\in\mathbb{R}} \bigl| P(F\leqslant x) -P(G
\leqslant x)\bigr|.
\]

\item The \textit{total variation distance} is
\[
d_{\mathrm{TV}}(F,G)= \sup_{A\in\mathscr{B}(\mathbb{R}) } \bigl| P(F\in A) -P(G\in A)\bigr|.
\]

\item The \textit{Fortet--Mourier distance} is
\[
d_{\mathrm{FM}}(F,G)= \sup_{\varphi\in\operatorname{Lip}(1), \|\varphi\|_\infty
\leqslant1} \bigl|E\bigl[ \varphi(F)\bigr]
-E\bigl[\varphi(G)\bigr]\bigr|.
\]
\end{itemize}

Plainly, $d_{\mathrm{W}}\geqslant d_{\mathrm{FM}}$ and $d_{\mathrm{TV}}\geqslant d_{\mathrm{Kol}}$. We recall
that the
topologies induced by $d_{\mathrm{W}}$, $d_{\mathrm{Kol}}$ and $d_{\mathrm{TV}}$, over the class of
probability measures on the real line, are strictly stronger than the
topology of convergence in distribution, whereas $d_{\mathrm{FM}}$ metrizes
convergence in distribution (see, e.g., \cite{np-book}, Appendix C,  for
a review of these facts).

\section{Quantitative stable convergence in dimension one}\label{sdq1}

We start by focussing on stable limits for one-dimensional Skorohod
integrals of order one, that is, random variables having the form $F =
\delta(u)$, where $u\in\mathbb{D}^{1,2}(\EuFrak{H})$. As already
discussed, this framework permits some interesting simplifications that
are not available for higher order integrals and higher dimensions.
Notice that any random variable $F $ such that $E[F]=0$ and
$E[F^{2}]<\infty$
can be written as  $F=\delta(u)$ for some $u\in\operatorname{Dom}\delta$. For
example, we can take $u=-DL^{-1}F$, or in the context of the standard
Brownian motion, we can take $u$ an adapted and square integrable process.

\subsection{Explicit estimates for smooth distances and stable CLTs}

The following estimate measures the distance between a Skorohod
integral of order 1, and a (suitably regular) mixture of Gaussian
distributions. In order to deduce a stable convergence result in the
subsequent Corollary~\ref{cstableclt}, we also consider an element
$I_1(h)$ in the first chaos of the isonormal process $X$.

\begin{prop}
\label{prop1-1}
Let $F\in\mathbb{D}^{1,2}$ be such that $E[F]=0$.
Assume $%
F=\delta(u)$ for some $u\in\mathbb{D}^{1,2}(\EuFrak{H})$. Let $S \ge0$
be such that $S^2\in
\mathbb{D}^{1,2}$, and let $\eta\sim\mathcal{N}(0,1)$ indicate a
standard Gaussian
random variable independent of the underlying isonormal Gaussian
process $X$. Let $h\in\mathfrak{H}$. Assume that $\varphi\dvtx \mathbb
{R}\rightarrow\mathbb{R}$ is $C^{3}$ with $\Vert\varphi^{\prime\prime}\Vert_{\infty},\Vert
\varphi^{\prime\prime\prime}\Vert_{\infty}<\infty$.
Then
%
\begin{eqnarray}
&& \bigl|E\bigl[\varphi\bigl(F + I_1(h)\bigr)\bigr] - E\bigl[
\varphi\bigl(S\eta+ I_1(h) \bigr)\bigr] \bigr|
\nonumber
\\
&&\label{smooth} \qquad \leqslant
\tfrac{1}{2}\bigl\Vert \varphi^{\prime\prime}\bigr\Vert_{\infty} E \bigl[2 \bigl|
\langle u, h\rangle_\mathfrak{H}\bigr|+\bigl|\langle u,DF\rangle_{\EuFrak{H}
}-S^{2}\bigr|\bigr]
\\
&&\qquad \quad
{}+\tfrac{1}{3}\bigl\Vert\varphi ^{\prime\prime\prime}\bigr\Vert _{\infty} E \bigl[\bigl|
\bigl\langle u,DS^{2}\bigr\rangle_{\EuFrak{H}}\bigr| \bigr].
\nonumber
\end{eqnarray}
\end{prop}

\begin{pf}
We proceed by interpolation. Fix $\varepsilon>0$
and set
$S_\varepsilon=\sqrt{S^2+\varepsilon}$. Clearly, $S_\varepsilon\in\mathbb
{D}^{1,2}$.
Let $g(t)=E[\varphi(I_1(h)+ \sqrt{t}F+%
\sqrt{1-t}S_\varepsilon\eta)]$, $t\in[0,1]$, and observe that
$E[\varphi
(F+ I_1(h))]-E[\varphi(S_\varepsilon%
\eta+ I_1(h))] = g(1)-g(0)=\int_0^1 g^{\prime}(t)\,dt$. For $t\in(0,1)$,
integrating by parts yields
\begin{eqnarray*}
g^{\prime}(t)&=&\frac{1}2 E \biggl[\varphi^{\prime}
\bigl(I_1(h)+ \sqrt {t}F+\sqrt {1-t}%
S_\varepsilon\eta
\bigr) \biggl( \frac{F}{\sqrt{t}}-\frac{S_\varepsilon\eta
}{\sqrt
{1-t}} \biggr) \biggr]
\\
&=&\frac{1}2 E \biggl[\varphi^{\prime}\bigl( I_1(h)+
\sqrt{t}F+\sqrt {1-t}S_\varepsilon\eta\bigr) \biggl( \frac{%
\delta(u)}{\sqrt{t}}-
\frac{S_\varepsilon\eta}{\sqrt{1-t}} \biggr) \biggr]
\\
&=&\frac{1}2 E \biggl[\varphi^{\prime\prime}\bigl( I_1(h)+
\sqrt{t}F+\sqrt {1-t}S_\varepsilon\eta\bigr)\\
&& \hspace*{21pt}{}\times \biggl(\frac{1}{\sqrt{t}}\langle
u, h\rangle _\mathfrak{H}+ \langle u,DF\rangle_\EuFrak H+
\frac{\sqrt{1-t}}{\sqrt{t}}\eta \langle u,DS_\varepsilon\rangle_\EuFrak H-
S_\varepsilon^2 \biggr) \biggr].
\end{eqnarray*}
Integrating again by parts with respect to the law of $\eta$ yields
\begin{eqnarray*}
g^{\prime}(t) &=&\frac{1}2 E \bigl[\varphi^{\prime\prime}
\bigl(I_1(h)+\sqrt{t}F+\sqrt {1-t}S_\varepsilon\eta\bigr)
\bigl(t^{-1/2} \langle u, h\rangle_\mathfrak{H}+ \langle u,DF
\rangle_\EuFrak H- S_\varepsilon^2 \bigr) \bigr]
\\
&& {}+\frac{1-t}{4\sqrt{t}} E \bigl[%
\varphi^{\prime\prime\prime}
\bigl(I_1(h)+\sqrt{t}F +\sqrt {1-t}S_\varepsilon\eta \bigr)\bigl
\langle u,DS ^2\bigr\rangle_\EuFrak H \bigr],
\end{eqnarray*}
where we have used the fact that $S_\varepsilon DS_\varepsilon= \frac{1}2
DS^2_\varepsilon=\frac{1}2 DS^2$.
Therefore,
\begin{eqnarray*}
&& \bigl\llvert E\bigl[\varphi\bigl(I_1(h)+F\bigr)\bigr]-E\bigl[\varphi
\bigl(I_1(h)+S_\varepsilon
\eta\bigr)\bigr] \bigr\rrvert\\
&& \qquad
\leqslant \frac{1}2 \bigl\|\varphi^{\prime\prime}\bigr\|_\infty E \bigl[2
\bigl| \langle u, h\rangle_\mathfrak{H}\bigr|+\bigl|\langle u,DF\rangle_\EuFrak H-
S ^2-\varepsilon\bigr| \bigr]
\\
&&\qquad \quad {}+\bigl\|\varphi^{\prime\prime\prime}\bigr\|_\infty E \bigl[%
\bigl|\bigl\langle
u,DS^2\bigr\rangle_\EuFrak H\bigr| \bigr]\int_0^1
\frac{1-t}{4\sqrt{t}}\,dt,
\end{eqnarray*}
and the conclusion follows letting $\varepsilon$ go to zero, because
$\int_0^1\frac{1-t}{4\sqrt{t}}\,dt=\frac{1}3$.

The following statement provides a stable limit theorem based on
Proposition~\ref{prop1-1}.

\begin{cor}\label{cstableclt} Let $S$ and $\eta$ be as in the
statement of Proposition~\ref{prop1-1}. Let $\{F_n\}$ be a sequence\vspace*{1pt} of
random variables such that $E[F_n] = 0$ and $F_n = \delta(u_n)$, where
$u_n\in\mathbb{D}^{1,2}(\EuFrak{H})$. Assume that the following conditions
hold as $n\to\infty$:
\begin{longlist}[3.]

\item[1.] $\langle u_n, DF_n\rangle_\EuFrak H\to S^2$ in $L^1(\Omega)$;

\item[2.] $\langle u_n, h\rangle_\EuFrak H\to0$ in $L^1(\Omega)$, for every
$h\in\EuFrak H$;

\item[3.] $\langle u_n, DS^2 \rangle_\EuFrak H\to0$ in $L^1(\Omega)$.
\end{longlist}
Then $F_n\stackrel{\mathrm{st}}{\rightarrow}S\eta$, and selecting $h=0$ in
(\ref{smooth}) provides
an upper bound for the rate of convergence of the difference $
|E[\varphi(F_n)] - E[\varphi(S\eta)] |$, for every $\varphi$
of class $C^3$ with bounded second and third derivatives.
\end{cor}

\begin{pf}
Relation (\ref{smooth}) implies that, if conditions 1--3 in the
statement hold true, then $ |E[\varphi(F_n + I_1(h))] -
E[\varphi(S\eta+ I_1(h) )] |\to0$ for every $h\in\EuFrak H$ and
every smooth test function $\varphi$. Selecting $\varphi$ to be a
complex exponential and using point 2 of Lemma~\ref{lstable} yields
the desired conclusion.
\end{pf}

\begin{rem}
(a) Corollary~\ref{cstableclt} should be compared with
Theorem~\ref{tHN} in the case $d=q=1$ (which exactly corresponds to
\cite{nounua}, Corollary~3.3). This result states that, if (i) $u_n\in
\mathbb{D}^{2,2}(\EuFrak{H})$ and (ii) $\{F_n\}$ is bounded in
$L^1(\Omega)$,
then it is sufficient to check conditions 1--2 in the statement of
Corollary~\ref{cstableclt} for some $S^2$ in $L^1(\Omega)$ in order to
deduce the stable convergence of $F_n$ to $S\eta$. The fact that
Corollary~\ref{cstableclt} requires more regularity on $S^2$, as well
as the additional condition 3, is compensated by the less stringent
assumptions on $u_n$, as well as by the fact that we obtain explicit
rates of convergence for a large class of smooth functions.

(b) The statement of \cite{nounua}, Corollary~3.3,  allows one
also to recover a modification of the so-called \textit{asymptotic Knight
Theorem} for Brownian martingales, as stated in \cite{RY}, Theorem~VIII.2.3. To see this, assume that $X$ is the isonormal Gaussian
process associated with a standard Brownian motion $B= \{B_t \dvtx
t\geqslant0\}$ [corresponding to the case $\EuFrak{H}= L^2(\mathbb{R}_+, ds)$] and
also that the
sequence $\{u_n\dvtx  n\geqslant1\}$ is composed of square-integrable processes
adapted to the natural filtration of $B$. Then, $F_n = \delta(u_n) =
\int_0^\infty u_n(s)\, dB_s$, where the stochastic integral is in the
It\^o sense, and the aforementioned asymptotic Knight theorem\vspace*{-1pt} yields that
the stable convergence of $F_n$ to $S\eta$ is implied by the following:
(A) $\int_0^t u_n(s)\, ds \stackrel{P}{\to} 0$, uniformly\vspace*{1pt} in $t$ in
compact sets and (B) $\int_0^\infty u_n(s)^2 \,ds \to S^2$ in
$L^1(\Omega)$.
\end{rem}

\subsection{Wasserstein and Kolmogorov distances}

The following statement provides a way to deduce rates of convergence
in the Wasserstein and Kolmogorov distance from the previous results.

\begin{teo}\label{tkol1}
Let $F\in\mathbb{D}^{1,2}$ be such that $E[F]=0$. Write $F=\delta
(u)$ for
some $u\in\mathbb{D}^{1,2}(\EuFrak{H})$. Let $S\in\mathbb{D}^{1,4}$, and
let $\eta\sim\mathcal{N}(0,1)$ indicate a standard Gaussian random
variable independent of the isonormal process $X$. Set
%
\begin{eqnarray}
\Delta&=&3 \biggl( \frac{1}{\sqrt{2\pi}}E \bigl[\bigl|\langle u,DF\rangle
_{\EuFrak H%
}-S^{2}\bigr| \bigr]+\frac{\sqrt{2}}{3} E \bigl[\bigl|\bigl
\langle u,DS^{2}\bigr\rangle _{\EuFrak H
}\bigr| \bigr]%
 \biggr)^{{1}/{3}}
\nonumber
\\
&&\label{delta}{}\times\max \biggl\{ \frac{1}{\sqrt{2\pi}}E \bigl[\bigl|\langle u,DF\rangle
_{\EuFrak H}-S^{2}\bigr| \bigr]+\frac{\sqrt{2}}{3} E \bigl[\bigl|\bigl
\langle u,DS^{2}\bigr\rangle _{\EuFrak H}\bigr|
 \bigr],\\
 && \hspace*{120pt}\nonumber\qquad \sqrt{
\frac{2}{\pi}} { \bigl( 2+E[S] \bigr)+E\bigl[|F|\bigr]} \biggr\} ^{{2}/{3}}.
\end{eqnarray}
Then $d_{\mathrm{W}}(F,S\eta)\leqslant\Delta$. Moreover, if there exists
$\alpha
\in(0,1]$ such that $E[|S|^{-\alpha}]<\infty$, then
%
\begin{equation}\label{kolmo}
d_{\mathrm{Kol}}(F,S\eta)\leqslant\Delta^{{\alpha}/({\alpha+1})} \bigl(1+E
\bigl[|S|^{-\alpha}\bigr] \bigr).
\end{equation}
\end{teo}

\begin{rem}
Theorem~\ref{tkol1} is specifically relevant whenever
one deals with sequences of random variables living in a finite sum of
Wiener chaoses. Indeed, in~\cite{nou-poly}, Theorem~3.1,  the following
fact is proved: let $\{F_n \dvtx  n\geqslant1\}$ be a sequence of random
variables living in the subspace $\bigoplus_{k=0}^p \mathcal{H}_k$, and
assume that $F_n$ converges in distribution to a nonzero random
variable $F_\infty$; then, there exists a finite constant $c>0$
(independent of $n$) such that
%
\begin{eqnarray}
&&  d_{\mathrm{TV}}(F_n, F_\infty)\leqslant c
d_{\mathrm{FM}}(F_n, F_\infty)^{{1}/({1+2p})} \leqslant c
d_{\mathrm{W}}(F_n, F_\infty)^{{1}/({1+2p})},
\nonumber
\\[-8pt]
\label{eig}\\[-8pt]
\eqntext{n \geqslant 1.}
\end{eqnarray}
Exploiting this estimate, and in the framework of random variables with
a finite chaotic expansion, the bounds in the Wasserstein distance
obtained in Theorem~\ref{tkol1} can be used to deduce rates of
convergence in total variation toward mixtures of Gaussian
distributions. The forthcoming Section~\ref{ssex1} provides an
explicit demonstration of this strategy, as applied to quadratic
functionals of a (fractional) Brownian motion.
\end{rem}

\begin{pf*}{Proof of Theorem~\ref{tkol1}}
It is divided into
two steps.

\begin{longlist}[Step 2.]
\item[\textit{Step} 1: \textit{Wasserstein distance}.] Let $\varphi\dvtx \mathbb
{R}\to\mathbb{R}$
be a function of class $C^3$ which is bounded together with all its
first three derivatives. For any $t\in(0,1)$, define
\[
\varphi_t(x)=\int_\mathbb{R}\varphi(\sqrt{t}y+
\sqrt{1-t}x)\,d\gamma(y),
\]
where $d\gamma(y)=\frac{1}{\sqrt{2\pi}}e^{-y^2/2}\,dy$ denotes the standard
Gaussian measure. Then, we may differentiate and integrate by parts to get
\begin{eqnarray*}
\varphi^{\prime\prime}_t(x) &=& \frac{1-t}{\sqrt{t}}\int
_\mathbb {R}y\varphi^{\prime
}(%
\sqrt{t}y+
\sqrt{1-t}x)\,d\gamma(y) \\
&=& \frac{1-t}{t}\int_\mathbb {R}
\bigl(y^2-1\bigr)\varphi (\sqrt{t%
}y+\sqrt{1-t}x)\,d\gamma(y)
\end{eqnarray*}
and
\[
\varphi^{\prime\prime\prime}_t(x)= \frac{(1-t)^{3/2}}{t}\int
_\mathbb{R} \bigl(y^2-1\bigr)\varphi^{\prime}(
\sqrt{t}y+\sqrt{1-t}x)\,d\gamma(y).
\]
Hence, for $0<t<1$ we may bound
%
\begin{equation}
\label{bound1}
\bigl\|\varphi^{\prime\prime}_t\bigr\|_\infty
\leqslant\frac{1-t}{\sqrt{t}}%
\bigl\|\varphi^{\prime}\bigr\|_\infty
\int_\mathbb{R}|y|\,d\gamma(y) \leqslant \sqrt{
\frac{2} \pi
} \frac{\|\varphi^{\prime}\|_\infty}{t}
\end{equation}
and
%
\begin{eqnarray}
\bigl\|\varphi^{\prime\prime\prime}_t\bigr\|_\infty
 & \leqslant & \frac
{(1-t)^{3/2}}{t}%
\bigl\|\varphi^{\prime}\bigr\|_\infty
\int_\mathbb{R}\bigl|y^2-1\bigr|\,d\gamma(y)
\nonumber
\\[-8pt]
\label{bound2}\\[-8pt]
\nonumber
& \leqslant &
\frac{%
\|\varphi^{\prime}\|_\infty}{t} \sqrt{\int_\mathbb {R}\bigl(y^2-1
\bigr)^2\,d\gamma(y)}=\frac
{%
\sqrt{2}\|\varphi^{\prime}\|_\infty}{t}.
\end{eqnarray}
Taylor expansion gives that
\begin{eqnarray*}
\bigl\vert E \bigl[\varphi(F)\bigr] - E \bigl[\varphi_t(F)\bigr]\bigr\vert&
\leqslant& \int_\mathbb{R}E \bigl[%
\bigl\llvert
\varphi( \sqrt{t} y + \sqrt{1-t} F) - \varphi(\sqrt {1-t} F) \bigr\rrvert \bigr]\, d
\gamma(y)
\\
&& {}+ E \bigl[\bigl\llvert \varphi( \sqrt{1-t} F) - \varphi( F) \bigr\rrvert
 \bigr]
\\
& \leqslant& \bigl\| \varphi^{\prime}\bigr\|_\infty\sqrt{t} \int
_\mathbb {R}| y | \,d\gamma(y) + \bigl\| \varphi^{\prime}
\bigr\|_\infty|\sqrt{1-t} -1| E \bigl[|F|\bigr]
\\
&\leqslant& \sqrt{t}\bigl\|\varphi^{\prime}\bigr\|_\infty \biggl\{ \sqrt {
\frac{2}\pi} + E\bigl[|F|\bigr] \biggr\}.
\end{eqnarray*}
Here, we used that $ |\sqrt{1-t}-1 |=t/(\sqrt
{1-t}+1)\leqslant
\sqrt{t}$. Similarly,
\begin{eqnarray*}
\bigl\vert E\bigl[ \varphi(S\eta)\bigr] - E\bigl[ \varphi_t(S\eta)\bigr]
\bigr\vert &\leqslant &  \sqrt{t}\bigl\|\varphi^{\prime}\bigr\|_\infty
\biggl\{ \sqrt{\frac{2} \pi} + E\bigl[|S\eta|\bigr] \biggr\} \\
&=&  \sqrt{
\frac{2} \pi} \sqrt{t}\bigl\|\varphi^{\prime}\bigr\|_\infty \bigl\{ 1
+ E[S] \bigr\}.
\end{eqnarray*}
Using \eqref{smooth} with \eqref{bound1}--\eqref{bound2} together
with the
triangle inequality and the previous inequalities, we have
%
\begin{eqnarray}
&& \bigl\vert E\bigl[\varphi(F)\bigr] - E\bigl[ \varphi(S\eta)\bigr] \bigr\vert\nonumber\\
&& \label{bound}\qquad \leqslant
\sqrt {t}\bigl\| \varphi^{\prime}\bigr\|_\infty \biggl( \sqrt{
\frac{2} \pi} \bigl\{ 2+ E[S] \bigr\}+ E\bigl[|F|\bigr] \biggr)
\\
&&\qquad \quad {}+\frac{\|\varphi^{\prime}\|_\infty}{t} \biggl\{ \frac{1} {\sqrt
{2\pi}}E \bigl[\bigl|\langle
u,DF\rangle_\EuFrak{H}-S^2\bigr| \bigr] +\frac{\sqrt{2}}{3} E
\bigl[\bigl|\bigl\langle u,DS^2\bigr\rangle_\EuFrak{H}\bigr| \bigr]
\biggr\}.\nonumber
\end{eqnarray}
Set
\[
\Phi_1={ \sqrt{\frac{2} \pi} \bigl\{ 2+ E[S] \bigr\} +
E\bigl[|F|\bigr]}
\]
and
\[
\Phi_2 = \frac{1} {\sqrt{2\pi}}E \bigl[\bigl|\langle u,DF
\rangle_\EuFrak H-S^2\bigr| \bigr] +
\frac{\sqrt{2}}{3} E
\bigl[\bigl|\bigl\langle u,DS^2\bigr\rangle_\EuFrak{H}\bigr| \bigr].
\]
The function $t\mapsto\sqrt{t} \Phi_1 + \frac{1}t \Phi_2$ attains its
minimum at $t_0=  ( \frac{2\Phi_2}{\Phi_1}  )^{2/3}$. Then,
if $%
t_0\leqslant1$ we choose $t=t_0$ and if $t_0 >1$ we choose $t=1$. With these
choices we obtain
%
\begin{eqnarray}
&& \bigl\vert E\bigl[\varphi(F)\bigr] - E\bigl[ \varphi(S\eta)\bigr]
\bigr\vert \nonumber
\\[-8pt]
\label{almostwasser}\\[-8pt]
\nonumber
&& \qquad \le  \bigl\|\varphi^{\prime}\bigr\|_\infty\Phi_2^{1/3}
\max\bigl(\bigl( 2^{-2/3}+ 2^{1/3}\bigr) \Phi_1^{2/3},
3 \Phi_2^{2/3} \bigr)
\leqslant  \bigl\|\varphi^{\prime}\bigr\|
_\infty\Delta.
\end{eqnarray}
This inequality can be extended to all
Lispchitz functions $\varphi$,
and this immediately yields that
$d_{\mathrm{W}}(F,S\eta)\leqslant\Delta$.

\item[\textit{Step} 2: \textit{Kolmogorov distance}.] Fix $z\in\mathbb{R}$
and $h>0$.
Consider the function $\varphi_h\dvtx \mathbb{R}\to[0,1]$ defined by
\[
\varphi_h(x) = \cases{
1, & \quad \mbox{if $x \leqslant z$},
\cr
0, &\quad  \mbox{if $x\geqslant z+h$},
\cr
\mathrm{linear},  & \quad  \mbox{if $z\leqslant x\leqslant z+h$},}
\]
and observe that $\varphi_h$ is Lipschitz with $\|\varphi^{\prime
}_h\|_\infty= 1/h$. Using that $\mathbf{1}_{(-\infty,z]}\leqslant
\varphi_h\leqslant\mathbf{1}_{(-\infty,z+h]}$ as well as (\ref{almostwasser}), we get
\begin{eqnarray*}
&& P[F\leqslant z]-P[S\eta\leqslant z]\\
&& \qquad \leqslant E\bigl[\varphi
_h(F)\bigr]-E\bigl[\mathbf{1}_{(-\infty,z]}(S\eta)\bigr]
\\
&& \qquad = E\bigl[\varphi_h(F)\bigr]-E\bigl[\varphi_h(S
\eta)\bigr]+E\bigl[\varphi_h(S\eta )\bigr]-E\bigl[\mathbf{1}
_{(-\infty,z]}(S\eta)\bigr]
\\
&& \qquad \leqslant  \frac{\Delta}{h} + P[z\leqslant S\eta\leqslant z+h].
\end{eqnarray*}
On the other hand, we can write
\begin{eqnarray*}
&&P[z\leqslant S\eta\leqslant z+h]
\\
&& \qquad =\frac{1}{\sqrt{2\pi}}\int_{\mathbb{R}^2} e^{-{x^2}/2}
\mathbf{1}_{[z,z+h]}(sx)\,dP_S(s)\,dx
\\
&& \qquad =\frac{1}{\sqrt{2\pi}} \biggl(\int_{\mathbb{R}_+}\, dP_S(s)
\int_{z/s}^{(z+h)/s}e^{-{x^2}/2}\,dx +\int
_{\mathbb{R}_-} \,dP_S(s) \int_{(z+h)/s}^{z/s}e^{-{x^2}/2}\,dx
\biggr)
\\
&&\qquad  \leqslant  \frac{|h|^\alpha}{\sqrt{2\pi}}\int_\mathbb {R}|s|^{-\alpha}
\,dP_S(s) \biggl(\int_\mathbb{R}e^{-{x^2}/({2(1-\alpha)})}\,dx
\biggr)^{1-\alpha}
\\
&& \qquad \leqslant |h|^\alpha E\bigl[|S|^{-\alpha}\bigr],
\end{eqnarray*}
because $(\int_\mathbb{R}e^{-{x^2}/({2(1-\alpha)})}\,dx
)^{1-\alpha}
= (\sqrt{1-\alpha}\int_\mathbb{R}e^{-{y^2}/{2}}\,dy
)^{1-\alpha}%
\leqslant\sqrt{2\pi}$, so that
\[
P[F\leqslant z]-P[S\eta\leqslant z]\leqslant\frac{\Delta}{h} +
|h|^\alpha E\bigl[|S|^{-\alpha}\bigr].
\]
Hence, by choosing $h=\Delta^{{1}/({\alpha+1})}$, we get that
\[
P[F\leqslant z]-P[S\eta\leqslant z]\leqslant\Delta^{{\alpha}/({\alpha+1})}%
 \bigl(1+E\bigl[|S|^{-\alpha}\bigr] \bigr).
\]
We prove similarly that
\[
P[F\leqslant z]-P[S\eta\leqslant z]\geqslant-
\Delta^{{\alpha}/({\alpha+1})}
 \bigl(1+E\bigl[|S|^{-\alpha}\bigr] \bigr),
\]
so the proof of (\ref{kolmo}) is done.\quad\qed
\end{longlist}
\noqed\end{pf*}

\subsection{Quadratic functionals of Brownian motion and fractional
Brownian motion}\label{ssex1}

We will now apply the results of the previous sections to some
nonlinear functionals of a fractional Brownian motion with Hurst
parameter $H\geqslant\frac{1}2$. Recall that a fractional Brownian motion
(fBm) with Hurst parameter $H\in(0,1)$
is a centered Gaussian process $B=\{B_{t}\dvtx  t\geqslant0\}$ with covariance
function
\[
E(B_{s}B_{t})= \tfrac{1}{2} \bigl(
t^{2H}+s^{2H}-|t-s|^{2H} \bigr).
\]
Notice that for $H=\frac{1}2$ the process $B$ is a standard Brownian motion.
We denote by
$\mathcal{E}$ the set of step functions on $[0,\infty) $. Let
$\EuFrak
H$ be the
Hilbert space defined as the closure of $\mathcal{E}$ with respect to the
scalar product
\[
\langle{\mathbf{1}}_{[0,t]},{\mathbf{1}}_{[0,s]} \rangle
_{%
\EuFrak H}=E(B_{s}B_{t}).
\]
The mapping $\mathbf{1}_{[0,t]}\rightarrow B_{t}$ can be extended to a
linear isometry between the Hilbert space $\EuFrak H$ and the Gaussian space
spanned by $B$.
We denote this isometry by $\phi\rightarrow B(\phi)$. In
this way, $\{B(\phi)\dvtx  \phi\in\EuFrak H\}$ is an isonormal Gaussian process.
In the case, $H>\frac{1}2$, the space $\EuFrak H$ contains all
measurable functions $\varphi\dvtx  \mathbb{R}_+ \rightarrow\mathbb{R}$
such that
\[
\int_0^\infty\!\int_0^\infty
\bigl|\varphi(s)\bigr| \bigl|\varphi(t)\bigr| |t-s| ^{2H-2} \,ds\,dt <\infty,
\]
and in this case if $\varphi$ and $\phi$ are functions satisfying this
integrability condition,
%
\begin{equation}
\label{eq7}
\langle\varphi, \phi \rangle_{\EuFrak H} = H(2H-1) \int
_0^\infty\!\int_0^\infty
\varphi(s) \phi(t) |t-s| ^{2H-2} \,ds\,dt.
\end{equation}
Furthermore, $L^{{1}/H}([0,\infty))$ is continuously embedded into
$\EuFrak H$. In what follows, we shall write
%
\begin{equation}
\label{ech} c_H = \sqrt{H(2H-1)\Gamma(2H-1)},\qquad H>1/2,
\end{equation}
and also $c_{{1}/2} := \lim_{H\downarrow {1}/{2}}c_H=$
$\frac{1}{\sqrt{2}}$.

The following statement contains explicit estimates in total variation
for sequences of quadratic Brownian functionals converging to a mixture
of Gaussian distributions. It represents a significant refinement of
\cite{py1}, Proposition~2.1  and \cite{petamwii}, Proposition~18.

\begin{teo}\label{to} Let $\{B_t \dvtx  t\geqslant0\}$ be a fBm of Hurst index
$H\geqslant\frac{1}2$. For every $n\geqslant1$, define
\[
A_n := \frac{n^{1+H}}{2} \int_0^1
t^{n-1} \bigl(B_1^2 - B_t^2
\bigr)\, dt.
\]
As $n\rightarrow\infty$, the sequence $A_n$ converges stably to
$S\eta$, where $\eta$ is a random variable independent of $B$ with law
$\mathcal{N}(0,1)$ and $S= c_H |B_1|$. Moreover, there exists a~constant $k$ (independent of $n$) such that
\[
d_{\mathrm{TV}}(A_n, S\eta) \leqslant k n^{ - ({1-H})/{15}}, \qquad  n
\geqslant 1.
\]
\end{teo}

The proof of Theorem~\ref{to} is based on the forthcoming Proposition~\ref{prop1} and Proposition~\ref{prop2}, dealing with the stable
convergence of some auxiliary stochastic integrals, respectively in the
cases $H=1/2$ and $H>1/2$. Notice that,\vspace*{1pt} since $\lim_{H\downarrow{1}/2} c_H=c_{{1}/{2}} =\frac{1}{\sqrt{2}}$, the statement of
Proposition~\ref{prop1} can be regarded as the limit of the statement
of Proposition~\ref{prop2}, as $H \downarrow\frac{1}2$.

\begin{prop} \label{prop1}
Let $B=\{B_t \dvtx  t\ge0\}$ be a standard Brownian motion. Consider the
sequence of It\^o integrals
\[
F_n= \sqrt{n}\int_{0}^{1}t^{n}B_{t}\,dB_{t},
\qquad n\geqslant1.
\]
Then the sequence $F_n$ converges stably to $S\eta$ as $n\rightarrow
\infty$, where $\eta$ is
a random variable independent of $B$ with law $\mathcal{N}(0,1)$ and
$S= \frac{|B_1|}{ \sqrt{2}}$. Furthermore, we have the following
bounds for the Wasserstein and Kolmogorov distances
\[
d_{\mathrm{Kol}}(F_{n},S\eta)\leqslant C_{\gamma}n^{-\gamma},
\]
for any $\gamma<\frac{1}{12}$, where $C_\gamma$ is a constant
depending on $\gamma$, and
\[
d_{\mathrm{W}}(F_{n},S\eta)\leqslant C n^{-{1}/6 },
\]
where $C$ is a finite constant independent of $n$.
\end{prop}

\begin{pf}
Taking into account that the Skorohod integral coincides with the It\^o
integral, we can write $F_n= \delta(u_n)$, where
$u_{n}(t)=\sqrt{n}t^{n}B_{t}\mathbf{1}_{[0,1]}(t)$. In order to apply
Theorem~\ref{tkol1}, we need to estimate the quantities
$E ( \llvert  \langle u_{n},\break DF_{n}\rangle_{\EuFrak H}-S^{2}\rrvert
)$ and
$ E  (\llvert  \langle u_{n},DS^{2}\rangle_{\EuFrak H} \rrvert
 )$.
We recall that $\EuFrak H=L^2( \mathbb{R}_+,ds)$. For $s\in[0,1]$, we
can write
\[
D_{s}F_{n}=\sqrt{n}s^{n}B_{s}+
\sqrt{n}\int_{s}^{1}t^{n}\,dB_{t}.
\]
As a consequence,
\[
\langle u_{n},DF_{n}\rangle_{\EuFrak H}=n\int
_{0}^{1}s^{2n}B_{s}^{2}\,ds+n
 \int_{0}^{1}s^n B_s
\biggl( \int_{s}^{1}t^{n}\,dB_{t}
\biggr) \,ds.
\]
From the estimates,
\begin{eqnarray*}
E \biggl( \biggl\llvert n\int_{0}^{1}s^{2n}B_{s}^{2}\,ds-
\frac
{B_{1}^{2}}{2}\biggr\rrvert \biggr) &\leqslant&n\int_{0}^{1}s^{2n}E
\bigl( \bigl\llvert B_{s}^{2}-B_{1}^{2}
\bigr\rrvert \bigr)\, ds+\biggl\llvert \frac{n}{2n+1}-\frac{1}{2}\biggr
\rrvert
\\
&\leqslant&2n \int_0^1 s^{2n}
\sqrt{1-s} \,ds+ \frac{1}{2(2n+1)}
\\
&\leqslant& \frac{2n} {\sqrt{2n+1}} \sqrt{\int_0^1
s^{2n} (1-s) \,ds}+ \frac{1}{2(2n+1)}
\\
&\leqslant&\frac{1}{\sqrt{2n}} + \frac{1}{4n}
\end{eqnarray*}
and
\begin{eqnarray*}
nE \biggl( \biggl\llvert \int_{0}^{1}s^{n}B_{s}
\biggl( \int_{s}^{1}t^{n}\,dB_{t}
\biggr) \,ds\biggr\rrvert \biggr) &\leqslant&\frac{n}{\sqrt{2n+1}}\int
_{0}^{1}s^{n+{1}/{2}}%
\sqrt{1-s^{2n+1}}\,ds
\\
&\leqslant&\frac{n}{(n+{3}/{2})\sqrt{2n+1}}\leqslant\frac
{1}{\sqrt{2n}},
\end{eqnarray*}
we obtain
%
\begin{equation}
\label{eq3}
E \bigl( \bigl\llvert \langle u_{n},DF_{n}
\rangle_{\EuFrak H}-S^{2}\bigr\rrvert \bigr) \leqslant
\frac{\sqrt{2}} {\sqrt{n} } + \frac{1}{4n}.
\end{equation}
On the other hand,
%
\begin{equation}
\label{eq4}
\bigl\llvert \bigl\langle u_{n},DS^{2}\bigr
\rangle_{\EuFrak H} \bigr\rrvert =\sqrt{n} E \biggl( \biggl\llvert
B_1 \int_{0}^{1}s^{n}B_{s}\,ds
\biggr\rrvert \biggr) \leqslant\frac{\sqrt
{n}}{n+{3}/{2}}%
\leqslant
\frac{1}{\sqrt{n}}.
\end{equation}
Notice that
%
\begin{equation}
\label{eq5}
E\bigl(|F_n |\bigr) \leqslant\frac{\sqrt{n}}{\sqrt{2n+2}} \leqslant
\frac{1}{
\sqrt{2}}.
\end{equation}
Therefore, using (\ref{eq3}), (\ref{eq4}) and (\ref{eq5}) and with
the notation of Theorem~\ref{tkol1}, for any constant $C<C_0$, where
\[
C_0=3 \biggl( \frac{1}{\sqrt{2\pi}} \biggl( \sqrt{2} +
\frac{1}4%
 \biggr) +\frac{\sqrt{2}}{3} \biggr) ^{{1}/{3}}
\biggl( \sqrt {\frac{2}{%
\pi}} \biggl( 2+ \frac{1}{\sqrt{\pi}}+
\frac{1}{\sqrt{2}} \biggr) \biggr) ^{{2}/{3}},
\]
there exists $n_0$ such that
for all $n\geqslant n_0$ we have
$\Delta\leqslant
Cn^{-{1}/{6}}$.
Therefore, $d_{\mathrm{W}}(F_{n},S\eta)\leqslant Cn^{-{1}/{6}}$ for $n\ge
n_0$. Moreover, $%
E[|S|^{-\alpha}]<\infty$ for any $\alpha<1$, which implies that
\[
d_{\mathrm{Kol}}(F_{n},S\eta)\leqslant C_{\gamma}n^{-\gamma},
\]
for any $\gamma<\frac{1}{12}$.
This completes the proof of the proposition.
\end{pf}

As announced, the next result is an extension of Proposition~\ref{prop1} to the case of the fractional Brownian motion with Hurst
parameter $H>\frac{1}2$.

\begin{prop} \label{prop2}
Let $B=\{B_t \dvtx  t\ge0\}$ be fractional Brownian motion with Hurst
parameter $H>\frac{1}2$. Consider the sequence of random variables $F_n
=\delta(u_n)$, $n\geqslant1$, where
\[
u_n(t)= n^H t^{n}B_{t}
\mathbf{1}_{[0,1]} (t).
\]
Then, the sequence $F_n$ converges stably to $S\eta$ as $n\rightarrow
\infty$, where $\eta$ is
a random variable independent of $B$ with law $\mathcal{N}(0,1)$ and
$S= c_H |B_1| $. Furthermore, we have the following bounds for the
Wasserstein and Kolmogorov distances
\[
d_{\mathrm{Kol}}(F_{n},S\eta)\leqslant C_{\gamma,H }n^{-\gamma},
\]
for any $\gamma<\frac{1-H}{6}$, where $C_{\gamma,H}$ is a constant
depending on $\gamma$ and $H$, and
\[
d_{\mathrm{W}}(F_{n},S\eta)\leqslant C_H
n^{-({1-H})/3},
\]
where $C_H$ is a constant depending on $H$.
\end{prop}

\begin{pf}
Let us compute
\[
D_{s}F_{n}=n^{H}s^{n}B_{s}+n^{H}
\int_{s}^{1}t^{n}\,dB_{t}.
\]
As a consequence,
\[
\langle u_{n},DF_{n}\rangle_{\EuFrak H}=\Vert
u_{n}\Vert_{\EuFrak H
}^{2}+n^{H} \biggl\langle u_{n},\int_{\cdot}^{1}t^{n}\,dB_{t}
\biggr\rangle_{\EuFrak H}.
\]
As in the proof of Proposition~\ref{prop1}, we need to estimate the
following quantities:
\[
\varepsilon_{n}=E \bigl(\bigl\llvert \Vert u_{n}
\Vert_{\EuFrak{H}}^{2}-S^{2}\bigr\rrvert \bigr)
\]
and
\[
\delta_{n}=E \biggl(\biggl\llvert n^{H} \biggl\langle
u_{n},\int_{\cdot
}^{1}t^{n}\,dB_{t}
\biggr\rangle_{\EuFrak H%
}\biggr\rrvert \biggr).
\]
We have, using (\ref{eq7}),
\begin{eqnarray*}
\varepsilon_{n} &\leqslant&H(2H-1) E \biggl( \biggl\llvert
2{n^{2H}}\int_{0}^{1}\!\int
_{0}^{t}s^{n}t^{n}B_{s}B_{t}(t-s)^{2H-2}\,ds\,dt-
\Gamma (2H-1)B_{1}^{2}\biggr\rrvert \biggr)
\\
&\leqslant&H(2H-1)n^{2H}E \biggl( \biggl\llvert 2\int
_{0}^{1}\!\int_{0}^{t}s^{n}t^{n}
\bigl[B_{s}B_{t}-B_{1}^{2}
\bigr](t-s)^{2H-2}\,ds\,dt
\biggr\rrvert \biggr)
\\
&&{}+H(2H-1)\biggl\llvert 2n^{2H}\int_{0}^{1}
\!\int_{0}^{t}s^{n}t^{n}(t-s)^{2H-2}\,ds\,dt-
\Gamma (2H-1)\biggr\rrvert
\\
&=&a_{n}+b_{n}.
\end{eqnarray*}
We can write for any $s\leqslant t$
\begin{eqnarray*}
E \bigl( \bigl\llvert B_sB_t - B_1^2
\bigr\rrvert \bigr) &=&  E \bigl( \bigl\llvert B_sB_t
-B_s B_1 + B_s B_1-
B_1^2 \bigr\rrvert \bigr) \\
&\leqslant & (1-t)^H +
(1-s)^H \leqslant2(1-s )^H.
\end{eqnarray*}
Using this estimate, we get
\[
a_{n } \leqslant4H(2H-1)n^{2H} \int_{0}^{1}
\!\int_{0}^{t}s^{n}t^{n}(1-s)^H
(t-s)^{2H-2}\,ds\,dt.
\]
For any positive integers $n,m$ set
%
\begin{equation}
\label{eq8}
\quad \rho_{n,m}= \int_{0}^{1}
\!\int_{0}^{t}s^{n}t^{m}
(t-s)^{2H-2}\,ds\,dt =\frac{\Gamma(n+1)\Gamma(2H-1)}{\Gamma(n+2H) (n+m +2H)}.
\end{equation}
Then, by H\"older's inequality,
\begin{eqnarray*}
a_n& \leqslant& 4H(2H-1)n^{2H} \rho_{n,n}^{1-H}
\biggl( \int_0^1\! \int_{0}^{t}s^{n}t^{n}
(1-s) (t-s)^{2H-2}\,ds\,dt \biggr)^H
\\
&=& 4H(2H-1)n^{2H} \rho_{n,n}^{1-H} (
\rho_{n,n}- \rho_{n+1,n} )^H.
\end{eqnarray*}
Taking into account that
\[
\rho_{n,n}- \rho_{n+1,n}= \frac{\Gamma(n+1)(n(2H+1) +4H^2)} {\Gamma
(n+2H)(2n+H) (n+2H)(2n+1+2H)},
\]
and using Stirling's formula, we obtain that $\rho_{n,n}$ is less than
{or} equal to a constant times $n^{-2H}$ and
$\rho_{n,n}- \rho_{n+1,n}$ is less than or equal to a constant times
$n^{-2H-1}$. This implies that
$a_{n} \leqslant C_H n^{-H}$, for some constant $C_H$ depending on $H$.

For the term $b_n$, using (\ref{eq8}) we can write
\[
b_{n}=H(2H-1)\Gamma(2H-1)\biggl\llvert \frac{2n^{2H}\Gamma(n+1)}{\Gamma
(n+2H)(2n+2H)}-1\biggr\rrvert,
\]
which converges to zero, by Stirling's formula, at the rate $ n^{-1}$.

On the other hand,
%
\begin{eqnarray}
\qquad \delta_{n} &=&H(2H-1)n^{2H}E \biggl( \biggl\llvert \int
_{0}^{1}\!\int_{0}^{1}s^{n}B_{s}
\biggl( \int_{t}^{1}r^{n}\,dB_{r}
\biggr) |t-s|^{2H-2}\,ds\,dt\biggr\rrvert \biggr)
\nonumber
\\[-8pt]
\label{x2}
\\[-8pt]
\nonumber
&\leqslant& H(2H-1)n^{2H}\int_{0}^{1}\!\int_{0}^{1}s^{n+H} \biggl[ E \biggl(
\biggl\llvert \int_{t}^{1}r^{n}\,dB_{r}
\biggr\rrvert ^{2} \biggr) \biggr] ^{1/2}|t-s|^{2H-2}\,ds\,dt.\hspace*{-6pt}
\end{eqnarray}
We can write, using the fact that $L^{{1}/H}([0,\infty))$ is
continuously embedded into $\EuFrak H$,
%
\begin{equation}\label{x1}
E \biggl( \biggl\llvert \int_{t}^{1}r^{n}\,dB_{r}
\biggr\rrvert ^{2} \biggr) \leqslant C_{H} \biggl( \int
_{t}^{1}r^{{n}/{H}}\,dr \biggr)
^{2H}\leqslant\frac
{C_{H}}{ ( {n}/{H}+1 ) ^{2H}}.
\end{equation}
Substituting (\ref{x11}) into (\ref{x22}) we obtain {$\delta
_{n}\leqslant
C_Hn^{H-1}$}, for some constant $C_H$, depending on $H$.  Thus,
\[
E \bigl( \bigl\llvert \langle u_{n},DF_{n}
\rangle_{\EuFrak H}-S^{2}\bigr\rrvert \bigr) \leqslant
C_H{n^{H-1}}.
\]
Finally,
\begin{eqnarray*}
E \bigl( \bigl\llvert \bigl\langle u_{n},DS^{2}\bigr\rangle_{\EuFrak H} \bigr\rrvert \bigr) &=&n^{H} E \biggl( \biggl
\llvert \int_{0}^{1}\!\int_{0}^{1}s^{n}B_{s}|t-s|^{2H-2}\,ds\,dt
\biggr\rrvert \biggr)
\\
&\leqslant&n^{H} \biggl\llvert \int_{0}^{1}
\!\int_{0}^{1}s^{n+H}|t-s|^{2H-2}\,ds\,dt
\biggr\rrvert \leqslant C_Hn^{H-1}.
\end{eqnarray*}
Notice that in this case  $E ( \llvert  \langle
u_{n},DF_{n}\rangle
_{\EuFrak H%
}-S^{2}\rrvert   ) $ converges to zero faster than  $E (
\llvert
\langle u_{n},DS^{2}\rangle_{\EuFrak H} \rrvert   ) $. As a
consequence, $\Delta
\leqslant C_{H}n^{ ({H-1})/{3}}$, for some constant $C_{H}$ and we
conclude the proof using Theorem~\ref{tkol1}.
\end{pf}

\begin{pf*}{Proof of Theorem~\ref{to}}
Using It\^o's formula (in
its classical form for \mbox{$H=\frac{1}2$}, and in\vspace*{1pt} the form discussed,
e.g., in \cite{nualartbook}, pages 293--294,  for the case $H>\frac{1}2$)
yields that
\[
\tfrac{1}{2} \bigl(B_1^2 - B_t^2
\bigr) = \delta \bigl(B_{\cdot} {\mathbf{1}}_{[t,1]}(\cdot) \bigr)
+\tfrac{1}{2}\bigl(1-t^{2H}\bigr)
\]
[note that $\delta (B_{\cdot} {\mathbf{1}}_{[t,1]}(\cdot) )$
is a
classical It\^o integral in the case $H=\frac{1}2$]. Interchanging
deterministic and stochastic integration by means of a stochastic
Fubini theorem yields therefore that
\[
A_n = F_n + H \frac{n^H}{2H+n}.
\]
In view of Propositions \ref{prop1} and \ref{prop2}, this implies that
$A_n$ converges in distribution to $S\eta$. The crucial point is now
that each random variable $A_n$ belongs to the direct sum $\mathcal
{H}_0 \oplus\mathcal{H}_2$: it follows that one can exploit the
estimate \eqref{eig} in the case $p=2$ to deduce that there exists a
constant $c$ such that
\[
d_{\mathrm{TV}}(A_n, S\eta) \leqslant c d_{\mathrm{W}}(A_n,S
\eta)^{{1}/5} \leqslant c \bigl(d_{\mathrm{W}}(F_n,S\eta)
+ d_{\mathrm{W}}(A_n,F_n) \bigr)^{{1}/5},
\]
where we have applied the triangle inequality. Since (trivially)
$d_{\mathrm{W}}(A_n, F_n) \leqslant H \frac{n^H}{2H+n} < n^{H-1}$, we deduce the
desired conclusion by applying the estimates in the Wasserstein
distance stated in Propositions \ref{prop1} and \ref{prop2}.
\end{pf*}

\section{Further notation and a technical lemma}\label{sec4}

\subsection{A technical lemma}

The following technical lemma is needed in the subsequent sections.

\begin{lemme}\label{lmipp}
Let $\eta_1,\ldots,\eta_d$ be a collection of i.i.d. $\mathcal{N}(0,1)$
random variables.
Fix $\alpha_1,\ldots,\alpha_d\in\mathbb{R}$ and integers
$k_1,\ldots,k_d\geqslant0$.
Then, for every $f \dvtx  \mathbb{R}^d \to\mathbb{R}$ of
class $C^{(k,\ldots,k)}$ (where $k = k_1+\cdots+k_d)$ such that $f$ and
all its partial derivatives have polynomial growth,
\begin{eqnarray*}
&& E \bigl[f(\alpha_1 \eta_1,\ldots,
\alpha_d\eta_d)\eta_1 ^{k_1}\cdots
\eta _d ^{k_d} \bigr]
\\
&& \qquad =\sum_{j_1=0}^{\lfloor k_1/2\rfloor}\cdots\sum
_{j_d=0}^{\lfloor
k_d/2\rfloor}\prod_{l=1}^d
\biggl\{\frac{k_l!}{%
2^{j_l}(k-2j_l)!j!}\alpha^{k_l-2j_l} \biggr\}
\\
&& \hspace*{62pt}\qquad \quad {}\times E \biggl[\frac{\partial
^{k_1+\cdots+k_d - 2(j_1+\cdots+j_d)}}{\partial x_1^{k_1-2j_1}\cdots
\partial x_d^{k_d-2j_d}}f(\alpha_1 \eta_1,
\ldots,\alpha_d\eta_d ) \biggr].
\end{eqnarray*}
\end{lemme}

\begin{pf}
By independence and conditioning, it suffices to
prove the claim for $d=1$, and in this case we write $\eta_1 = \eta$,
$k_1 = k$, and so on. The decomposition of the random variable $\eta
^{k}$ in terms of Hermite polynomials is given by
\[
\eta^{k}=\sum_{j=0}^{\lfloor k/2\rfloor}
\frac{k!}{2^{j}(k-2j)!j!}
H_{k-2j}(\eta),
\]
where $H_{k-2j}(x)$ is the $%
(k-2j)$th Hermite polynomial. Using the relation $E[f(\alpha\eta
)H_{k-2j}(\eta)] = \alpha^{k-2j}E[f^{(k-2j)}(\alpha\eta)]$, we deduce
the desired conclusion.
\end{pf}

\subsection{Notation}\label{snot}

The following notation is needed in order to state our next results.
For the rest of this section, we fix integers $m\geqslant0$ and
$d\geqslant1$.
\begin{longlist}[(iii)]

\item[(i)] In what follows, we shall consider smooth functions
%
\begin{equation}
\label{esf}
\qquad\psi\dvtx  \mathbb{R}^{m\times d} \to\mathbb{R}\dvtx  (y_1,
\ldots,y_m; x_1,\ldots,x_d) \mapsto\psi
(y_1,\ldots,y_m ; x_1,\ldots,x_d).
\end{equation}
Here, the implicit convention is that, if $m=0$, then $\psi$ does not
depend on $(y_1,\ldots,y_m)$. We also write
\[
\psi_{x_k} = \frac{\partial}{\partial x_k}\psi, \qquad  k=1,\ldots,d.
\]

\item[(ii)] For every integer $q\geqslant1$, we write $\mathscr{A}(q)
=\mathscr{A}(q; m,d)$ (the dependence on $m,d$ is dropped whenever
there is no risk of confusion) to indicate the collection of all
$(m+q(1+d))$-dimensional vectors with nonnegative integer entries of
the type
%
\begin{equation}
\label{eaq}
\alpha^{(q)} = (k_1,\ldots,k_q;
a_1,\ldots,a_m ; b_{ij}, i=1,\ldots,q, j=1,
\ldots,d),
\end{equation}
verifying the set of Diophantine equations
%
\begin{equation}
\label{Dio}
%
\begin{array} {rll} k_1+2k_2+\cdots+
qk_q &=&q,
\\
a_1+\cdots+a_m + b_{11}+\cdots+b_{1d}
&=&k_1,
\\
b_{21}+\cdots+b_{2d} &=&k_2,
\\
\cdots&&
\\
b_{q1}+\cdots+b_{qd}&=&k_q.
\end{array}
%
\end{equation}

\item[(iii)] Given $q\geqslant1$ and $\alpha^{(q)}$ as in (\ref{eaq}), we define
%
\begin{equation}
\label{ecaq}
C\bigl(\alpha^{(q)}\bigr) := \frac{q!}{\prod_{i=1}^q i!^{k_i}\prod_{l=1}^m a_l!
\prod_{i=1}^q\prod_{j=1}^db_{ij}!}.
\end{equation}

\item[(iv)] Given a smooth function $\psi$ as in (\ref{esf}) and a
vector $\alpha^{(q)} \in\mathscr{A}(q)$ as in (\ref{eaq}), we set
%
\begin{equation}
\label{edaq}
\partial^{\alpha^{(q)}} \psi:= \frac{\partial^{k_1+\cdots+
k_d}}{\partial y_1^{a_1}\cdots\partial y_m^{a_m}\partial x_1^{b_{11}
+\cdots+ b_{q1}} \cdots\partial x_d^{b_{1d} +\cdots+ b_{qd}}} \psi.
\end{equation}
The coefficients $C(\alpha^{(q)})$ and the differential operators
$\partial^{\alpha^{(q)}}$, defined respectively in (\ref{ecaq}) and
(\ref{edaq}), enter the generalized Faa di Bruno formula (as proved,
e.g., in~\cite{Mis}) that we will use in the proof of our main results.

\item[(v)] For every integer $q\geqslant1$, the symbol $\mathscr{B}(q)
=\mathscr{B}(q; m,d)$ indicates the class of all
$(m+q(1+2d))$-dimensional vectors with nonnegative integer entries of
the type
%
\begin{equation}
\label{ebq}
\beta^{(q)} = \bigl(k_1,
\ldots,k_q; a_1,\ldots,a_m;
b'_{ij}, b''_{ij},
i=1,\ldots,q, j=1,\ldots,d\bigr),
\end{equation}
such that
%
\begin{equation}
\label{eabq}
\qquad\alpha\bigl(\beta^{(q)}\bigr) := \bigl(k_1,
\ldots,k_q ; a_1,\ldots,a_m ;
b'_{ij} + b''_{ij},
i=1,\ldots,q, j=1,\ldots,d\bigr),
\end{equation}
is an element of $\mathscr{A}(q)$, as defined at point (ii). Given
$\beta^{(q)}$ as in (\ref{ebq}), we also adopt the notation
%
\begin{eqnarray}
\bigl| b' \bigr| &: =&  \sum_{i=1}^q
\sum_{j=1}^d b'_{ij},\qquad
\bigl| b'' \bigr| : = \sum_{i=1}^q
\sum_{j=1}^d b''_{ij},
\nonumber
\\[-8pt]
\label{enorms}\\[-8pt]
\nonumber
\bigl| b''_{\bullet j} \bigr|  &: =&  \sum
_{i=1}^q b''_{ij},
\qquad j=1,\ldots,d.
\end{eqnarray}

\item[(vi)] For every $\beta^{(q)}\in\mathscr{B}(q)$ as in (\ref{ebq})
and every $(l_1,\ldots,l_d)$ such that $l_s \in\{0,\ldots,\lfloor|
b''_{\bullet s}|/2\rfloor\}$, $s=1,\ldots,d$, we set
%
\begin{eqnarray}
&& W\bigl(\beta^{(q)}; l_1,
\ldots,l_d\bigr)
\nonumber
\\[-8pt]
\label{eW}\\[-8pt]
\nonumber
&& \qquad := C\bigl(\alpha\bigl(\beta^{(q)}\bigr)
\bigr)\prod_{i=1}^q\prod
_{j=1}^d \pmatrix{b'_{ij}+b''_{ij}
\cr
b'_{ij}}\prod_{s=1}^d
\frac{|
b''_{\bullet s}| !}{2^{l_s}(| b''_{\bullet s}| - 2l_s)!l_s!},
\end{eqnarray}
where $C(\alpha(\beta^{(q)}))$ is defined in (\ref{ecaq}), and
%
\begin{equation}
\label{edstar}
\partial_{\star}^{(\beta^{(q)};l_1,\ldots,l_d)} :=\partial^{\alpha
(\beta
^{(q)})}
\frac{{\partial} ^ { |b''| - 2(l_1+\cdots+l_d) }}{\partial
x_1^{|b''_{\bullet1}| - 2l_1}\cdots\partial x_d^{|b''_{\bullet d}| - 2l_d}},
\end{equation}
where $\alpha(\beta^{(q)})$ is given in (\ref{eabq}), and $\partial
^{\alpha(\beta^{(q)})}$ is defined according to (\ref{edaq}).

\item[(vii)] The Beta function $B(u,v)$ is defined as
\[
B(u,v) = \int_0^1 t^{u-1}(1-t)^{v-1}
\,dt, \qquad u,v>0.\
\]
\end{longlist}

\section{Bounds for general orders and dimensions}\label{sdqany}

\subsection{A general statement}

The following statement contains a general upper bound, yielding stable
limit theorems and associated explicit rates of convergence on the
Wiener space.

\begin{teo}\label{tmain}
Fix integers $m\geqslant0$, $d\geqslant1$ and
$q_j\geqslant1$, $j=1,\ldots,d$. Let $\eta= (\eta_1,\ldots,\eta_d)$ be a vector
of i.i.d. $\mathscr{N}(0,1)$ random variables independent of the
isonormal Gaussian process $X$. Define $\hat{q} = \max_{j=1,..,d} q_j$.
For every $j=1,\ldots,d$, consider a symmetric random element $u_j\in
\mathbb{D}^{2\hat{q},4\hat{q}}(\mathfrak{H}^{2q_j})$, and introduce the
following notation:
\begin{itemize}
\item $F_j := \delta^{q_j}(u_j)$ and $F := (F_1,\ldots,F_d)$;

\item $(S_1,\ldots,S_d)$ is a vector of real-valued elements of
$\mathbb{D}^{\hat{q},4\hat{q}}$, and
\[
S\cdot\eta:= (S_1\eta_1,\ldots,S_d
\eta_d).
\]
\end{itemize}
Assume that the function $\varphi\dvtx \mathbb{R}^{m\times d} \rightarrow
\mathbb{R}$
admits continuous and bounded partial derivatives up to the order
$2\hat
{q}+1$. Then, for every $h_1,\ldots,h_m\in\EuFrak H$,
%
\begin{eqnarray}
&& \bigl\llvert E\bigl[\varphi\bigl(X(h_1),\ldots,X(h_m);
F\bigr)\bigr] -E\bigl[\varphi \bigl(X(h_1),\ldots,X(h_m);
S\cdot\eta\bigr)\bigr]\bigr\rrvert
\nonumber
\\
&& \label{w1}\qquad\leqslant\frac{1}2 \sum_{k,j=1}^d
\biggl\llVert \frac{\partial
^2}{\partial x_k\,
\partial x_j} \varphi\biggr\rrVert _\infty E \bigl[
\bigl\llvert \bigl\langle D^{q_k}F_j , u_k
\bigr\rangle_{\EuFrak H^{\otimes q_k}} - {\mathbf{1}}_{j=k} S^2_j
\bigr\rrvert \bigr]
\\
&& \qquad\quad {}+ \frac{1}2 \sum_{k=1}^d \sum
_{\beta^{(q_k)}\in\mathscr
{B}_0(q_k)} \sum_{l_1=0}^{\lfloor|b''_{\bullet1}|/2\rfloor}
\cdots \sum_{l_d=0}^{\lfloor|b''_{\bullet d}|/2\rfloor} \widehat{W}\bigl(
\beta ^{(q_k)}; l_1,\ldots,l_d\bigr)
\nonumber
\\[-8pt]
\label{w2}\\[-8pt]
\nonumber
&& \quad \qquad {}\times\bigl\llVert
\partial_\star^{(\beta
^{(q_k)};l_1,\ldots,l_d)}\varphi_{x_k} \bigr\rrVert
_\infty
\\
&& \qquad\quad{}\times E \Biggl[\prod_{s=1}^d
S^{ |b''_{\bullet s}| -2l_s}\nonumber\\
&& \qquad\quad{}\times \Biggl\llvert \Biggl\langle u_k,
h_1^{\otimes a_1}\otimes\cdots
\otimes h_m^{\otimes a_m}
\bigotimes_{i=1}^{q_k} \bigotimes
_{j=1}^{d} \bigl\{ \bigl(D^iF_j
\bigr)^{\otimes b'_{ij}}\otimes\bigl(D^iS_j
\bigr)^{\otimes b''_{ij}} \bigr\} \Biggr\rangle_{\EuFrak H^{\otimes q_k}} \Biggr\rrvert
\Biggr],
\nonumber
\end{eqnarray}
where we have adopted the same notation as in Section~\ref{snot}, with
the following additional conventions: \textup{(a)} $\mathscr{B}_0(q)$ is the
subset of $\mathscr{B}(q)$ composed\vspace*{1.5pt} of those $\beta(q_k)$ as in (\ref{ebq}) such that $b'_{qj} = 0$ for $j=1,\ldots,d$, \textup{(b)} $\widehat
{W}(\beta
^{(q_k)}; l_1,\ldots,l_d) := W(\beta^{(q_k)}; l_1,\ldots,l_d)\times{B(|b'|/2
+1/2; |b''| /2+1)}$, where $B$ is the Beta\ function.
\end{teo}

\subsection{Case $m=0$, $d=1$}

Specializing Theorem~\ref{tmain} to the choice of parameters $m=0$,
$d=1$ and $q\geqslant1$ yields the following estimate on the distance
between the laws of a (multiple) Skorohod integral and of a mixture of
Gaussian distributions.

\begin{prop} \label{p1}
Suppose that $u\in\mathbb{D}^{2q,4q}(\mathfrak{H}^{2q})$ is symmetric.
Let $%
F=\delta^{q}(u)$. Let $S\in\mathbb{D}^{q,4q}$, and let $\eta\sim
\mathcal{N}(0,1)$
indicate a standard Gaussian random variable, independent of the underlying
isonormal process $X$. Assume that $\varphi\dvtx \mathbb{R}\rightarrow
\mathbb{R}$
is $C^{2q+1}$ with $\Vert\varphi^{(k)}\Vert_{\infty}<\infty$ for
any $%
k=0,\ldots,2q+1$. Then
\begin{eqnarray*}
&& \bigl|E\bigl[\varphi(F)\bigr]-E\bigl[\varphi(S\eta)\bigr] \bigr|\\
&&\qquad \leqslant
\frac{1}{2}\bigl\Vert\varphi^{\prime\prime}\bigr\Vert_{\infty} E
 \bigl[\bigl|{\bigl\langle u,D^qF\bigr\rangle_{\EuFrak H^{\otimes q}}}-S^{2}\bigr|
\bigr]
\\
&&\qquad \quad {}+\sum_{(b',b'')\in\mathcal{Q},b'_{q}=0} \sum_{j=0}^{ {\lfloor
\llvert  b'' \rrvert  /2 \rfloor%
} }
c_{q,b',b'',j} \bigl\llVert \varphi^{(1+\llvert  b'\rrvert  +2\llvert
b''\rrvert
-2j)}\bigr\rrVert
_{\infty}
\\
&& \qquad \qquad \hspace*{90pt} {}\times E \bigl[S^{ |b''|-2j}\\
&& \qquad \qquad \hspace*{90pt} \hspace*{23pt}{}\times \bigl\llvert \bigl\langle u, ( DF )
^{\otimes
b'_{1}}\otimes\cdots\otimes \bigl( D^{q-1}F \bigr)
^{\otimes
b'_{q-1}}\\
&& \hspace*{90pt}\hspace*{53pt}\qquad \qquad {}\otimes ( DS ) ^{\otimes b''_{1}}\otimes\cdots
\otimes \bigl(
D^{q}S \bigr) ^{\otimes b''_{q}} \bigr\rangle_{\mathfrak
{H}^{\otimes
q}} \bigr
\rrvert \bigr],
\end{eqnarray*}
where $\mathcal{Q}$ is the set of all pairs of $q$-ples
$b'=(b'_{1},b'_{2},\ldots,b'_{q})$ and $b''=(b''_{1},\ldots,b''_{q})$ of nonnegative integers
satisfying the constraint $b'_{1}+2b'_{2}+\cdots
+qb'_{q}+b''_{1}+2b''_{2}+\cdots
+qb''_{q}=q$.  The constants $c_{q,b',b'',j}$ are given by
\begin{eqnarray*}
 c_{q,b',b'',j}
&=& \frac{1}2 B\bigl(\bigl|b'\bigr|/2 +1/2,
\bigl|b''\bigr|/2+1\bigr) \\
&& {}\times\prod_{i=1}^q
{b_i \choose b'_i} \times
\frac{ |b''|!}{ 2^j (|b''| -2j)! j!} \times\frac{ q!}{ \prod_{i=1}^q
i!^{b_i} b_i!},
\end{eqnarray*}
where $b=b'+b''$.
\end{prop}

In the particular case $q=2$ we obtain the following result.

\begin{prop} \label{p11}
Suppose that $u\in\mathbb{D}^{4,8}(\mathfrak{H}^{4})$ is symmetric.
Let $%
F=\delta^{2}(u)$. Let $S\in\mathbb{D}^{2,8}$, and let $\eta\sim
\mathcal{N}(0,1)$
indicate a standard Gaussian random variable, independent of the underlying
isonormal process $X$. Assume that $\varphi\dvtx \mathbb{R}\rightarrow
\mathbb{R}$
is $C^{5}$ with $\Vert\varphi^{(k)}\Vert_{\infty}<\infty$ for any $
k=0,\ldots,5$. Then
\begin{eqnarray*}
&& \bigl|E\bigl[\varphi(F)\bigr]-E\bigl[\varphi(S\eta)\bigr] \bigr|  \\
&& \qquad \leqslant
\frac{1}{2}\bigl\Vert\varphi^{\prime\prime}\bigr\Vert_{\infty} E
 \bigl[\bigl|\bigl\langle u,D^2F\bigr\rangle_{{\EuFrak H}^\otimes2}-S^{2}\bigr|
\bigr]
\\
&& \quad \qquad{}+ {C_0} \max_{3\leqslant i \leqslant5} \bigl\|\varphi^{(i)}
\bigr\|_\infty \bigl( E \bigl[ \bigl\llvert \bigl\langle u,
(DF)^{\otimes2} \bigr\rangle _{{\EuFrak
H}^{\otimes2}} \bigr\rrvert \bigr] + E \bigl[
S \bigl\llvert \langle u, DF\otimes DS \rangle_{{\EuFrak H}^{\otimes2}} \bigr\rrvert
\bigr]
\\
&& \hspace*{91pt}\quad \qquad{}+ E \bigl[ \bigl(S^2+1\bigr) \bigl\llvert \bigl\langle u,
(DS)^{\otimes2} \bigr\rangle_{{\EuFrak H}^{\otimes2}} \bigr\rrvert \bigr]\\
 && \hspace*{120pt}\hspace*{91pt}\quad \qquad{}+ E \bigl[S
\bigl\llvert \bigl\langle u, D^2S \bigr\rangle_{{\EuFrak
H}^{\otimes2}} \bigr
\rrvert \bigr] \bigr),
\end{eqnarray*}
where $C_0= \frac{1}2 B(\frac{1}2, \frac{3}2) + \frac{3}2 B(\frac{3}2,1)+
B(\frac{1}2,2)$.
\end{prop}

Taking into account that $DS^2= 2SDS$ and $D^2S^2= 2DS\otimes DS+
2SD^2S$, we can write the above estimate in terms of the derivatives of
$S^2$, which is helpful in the applications. In this way, we obtain
\begin{eqnarray}
&&\bigl|E\bigl[\varphi(F)\bigr]-E\bigl[\varphi(S\eta)\bigr] \bigr|
\nonumber\\
&& \qquad \leqslant
\frac{1}{2}\bigl\Vert\varphi^{\prime\prime}\bigr\Vert_{\infty} E
 \bigl[\bigl|\bigl\langle u,D^2F\bigr\rangle_{{\EuFrak H}^{\otimes2}}-S^{2}\bigr|
\bigr]\nonumber
\\
&& \label{eq11}\quad\qquad {}+C_0 \max_{3\leqslant i \leqslant5}  \bigl\|\varphi^{(i)}
\bigr\| _\infty \bigl( E \bigl[ \bigl\llvert \bigl\langle u,
(DF)^{\otimes2} \bigr\rangle _{{\EuFrak
H}^{\otimes2}} \bigr\rrvert \bigr] + E \bigl[
\bigl\llvert \bigl\langle u, DF\otimes{DS^2} \bigr
\rangle_{{\EuFrak H}^{\otimes2}} \bigr\rrvert \bigr]
\\
&&\hspace*{91pt} \quad\qquad {}+ E \bigl[ \bigl(S^{-2}+1\bigr) \bigl\llvert \bigl\langle u, {
\bigl(DS^2\bigr)^{\otimes2}} \bigr\rangle_{{\EuFrak H}^{\otimes2}} \bigr
\rrvert \bigr] \nonumber\\
&&\hspace*{118pt}\hspace*{91pt} \quad\qquad {}+ E \bigl[ \bigl\llvert \bigl\langle u, {D^2S
^2} \bigr\rangle_{{\EuFrak
H}^{\otimes2}} \bigr\rrvert \bigr]
\bigr).\nonumber
\end{eqnarray}
Notice that a factor $S^{-2}$ appears in the right-hand side of the
above inequality.

\subsection{Case $m>0$, $d=1$}
Fix $q\geqslant1$. In the case $m>0$, $d=1$, the class $\mathscr
{B}(q)$ is
the collection of all vectors with nonnegative integer entries of the
type $\beta^{(q)} = (a_1,\ldots,a_m; b'_1,b''_1,\ldots,b'_q,b''_q)$ verifying
\[
a_1+\cdots+a_m +\bigl(b'_1+b''_1
\bigr)+\cdots+q\bigl(b'_q+b''_q
\bigr) = q,
\]
whereas $\mathscr{B}_0(q)$ is the subset of $\mathscr{B}(q)$ verifying
$b'_q = 0$. Specializing Theorem~\ref{tmain} yields upper bounds for
one-dimensional $\sigma(X)$-stable convergence.\vspace*{-2pt}

\begin{prop} \label{p2}
Suppose that $u\in\mathbb{D}^{2q,4q}(\mathfrak{H}^{2q})$ is symmetric,
select $h_1,\ldots, h_m\in\EuFrak H$, and write $\mathbf{X} = (X(h_1),\ldots,X(h_m))$.
Let $%
F=\delta^{q}(u)$. Let $S\in\mathbb{D}^{q,4q}$, and let $\eta\sim
\mathcal{N}(0,1)$
indicate a standard Gaussian random variable, independent of the underlying
Gaussian field $X$. Assume that
\[
\varphi\dvtx \mathbb{R}^m\times\mathbb{R}\rightarrow\mathbb{R}\dvtx
(y_1,\ldots,y_m,x) \mapsto\varphi(y_1,
\ldots,y_m,x)
\]
admits continuous and bounded partial derivatives up to the order
$2q+1$. Then
\begin{eqnarray*}
&& \bigl|E\bigl[\varphi(\mathbf{X}, F)\bigr]-E\bigl[\varphi(\mathbf{X}, S\eta)\bigr]
\bigr|
\\
&& \qquad \leqslant\frac{1}{2}\biggl\llVert \frac{\partial^2}{\partial x^2} \varphi\biggr
\rrVert _{\infty} E
 \bigl[\bigl|\bigl\langle u,D^qF\bigr
\rangle_{{\EuFrak H}^{\otimes q}}-S^{2}\bigr| \bigr] \\[2pt]
&& \qquad \quad {}+\frac{1}{2} \sum
_{{\beta^{(q)}}\in\mathscr{B}_0(q)} \sum_{j=0}^{\lfloor\llvert  b''\rrvert /2
\rfloor}
\widehat{W}\bigl(\beta^{(q)},j\bigr)
\biggl\llVert \frac{\partial^{|a|}}{\partial y_1^{a_1}\cdots
\partial y_m^{a_m}}\frac{\partial^{1+|b'|+2|b''|-2j}}{\partial
x^{1+|b'|+2|b''|-2j}} \varphi\biggr\rrVert
_{\infty}
\\[2pt]
&& \qquad \qquad{}\times E \Biggl[ S^{|b''|-2j} \Biggl\llvert \Biggl\langle u,
h_1^{\otimes
a_1}\otimes\cdots\\[2pt]
&& \qquad \hspace*{104pt}{}\otimes h_m^{\otimes a_m}
\bigotimes_{i=1}^{q} \bigl\{
\bigl(D^iF\bigr)^{\otimes b'_{i}}\otimes\bigl(D^iS
\bigr)^{\otimes b''_{i}} \bigr\} \Biggr\rangle_{\EuFrak H^{\otimes q}} \Biggr\rrvert
\Biggr],
\end{eqnarray*}
where $|a| = a_1+\cdots+a_m$.
\end{prop}

\subsection{Proof of Theorem \texorpdfstring{\protect\ref{tmain}}{5.1}}

The proof is based on the use of an interpolation argument. Write $\mathbf{X} = (X(h_1),\ldots,X(h_m))$ and $g(t)=E[\varphi(\mathbf{X}; \sqrt{t}F+
\sqrt{1-t} S\cdot\eta)]$, $t\in{}[0,1]$, and observe\vspace*{1pt} that
$E[\varphi
(\mathbf{X}; F)]-E[\varphi(\mathbf{X} ; S\eta)]=g(1)-g(0)=\int_{0}^{1}g^{\prime}(t)\,dt$. For $t\in
(0,1)$, by integrating by parts with respect either to $F$ or
to $\eta$, we get
\begin{eqnarray*}
g^{\prime}(t) &=&\frac{1}{2} \sum_{k=1}^d
E \biggl[ \varphi _{x_k}(\mathbf{X} ; \sqrt{t}F+\sqrt{1-t%
}S
\cdot\eta) \biggl( \frac{F_k}{\sqrt{t}}-\frac{S_k \eta_k }{\sqrt
{1-t}} \biggr) \biggr]
\\[2pt]
&=&\frac{1}{2} \sum_{k=1}^d E
\biggl[ \varphi_{x_k}(\mathbf{X}; \sqrt {t}F+\sqrt{1-t%
}S\cdot
\eta) \biggl( \frac{\delta^{q_k}(u_k)}{\sqrt{t}}-\frac{S_k
\eta
_k }{\sqrt{1-t}} \biggr) \biggr]
\\[2pt]
&=&\frac{1}{2\sqrt{t}} \sum_{k=1}^d E
\bigl[ \bigl\langle D^{q_k}\varphi_{x_k}(\mathbf{X};
\sqrt{t}%
F+\sqrt{1-t}S\cdot\eta),u_k \bigr
\rangle_{\mathfrak{H}^{\otimes
q_k}} \bigr]
\\[2pt]
&&{}-\frac{%
1}{2}\sum_{k=1}^d E
\biggl[\frac{\partial^2}{\partial x_k^2} \varphi (\mathbf{X}; \sqrt{t}F+\sqrt{1-t}S\cdot
\eta)S_k^{2}%
 \biggr].
\end{eqnarray*}
Using the Faa di Bruno formula for the iterated derivative of the
composition of a~function with a vector of functions (see \cite{Mis}, Theorem~2.1), we infer that, for every $k=1,\ldots,d$,
%
\begin{eqnarray}
&&\hspace*{6pt}\bigl\langle D^{q_k}\varphi_{x_k} (\mathbf{X} ; \sqrt{t}F+
\sqrt{1-t}S\cdot \eta) , u_k \bigr\rangle_{\EuFrak H^{\otimes q_k}}
\nonumber
\\[2pt]
&&\hspace*{6pt}\label{e1} \qquad=\sum_{\alpha^{(q_k)} \in\mathscr{A}(q_k) } C\bigl(\alpha^{(q_k)}\bigr)
\partial^{(\alpha^{(q_k)})}\varphi_{x_k} (\mathbf{X} ; \sqrt{t}F+\sqrt
{1-t}S\cdot\eta)
\\[2pt]
&&\hspace*{6pt} \quad\qquad {}\times \Biggl\langle h_1^{\otimes
a_1} \otimes\cdots\otimes
h_m^{\otimes a_m}\bigotimes_{i=1}^{q_k}
\bigotimes_{j=1}^d \bigl( D^i
(\sqrt{t} F_j + \sqrt{1-t} S_j\eta_j)
\bigr)^{\otimes b_{ij}} , u_k \Biggr\rangle_{\EuFrak H^{\otimes
q_k}}.
\nonumber
\end{eqnarray}\vspace*{-12pt}
\noindent For every $i=1,\ldots,q_k$, every $j=1,\ldots,d$ and every symmetric $v \in
\EuFrak H
^{\otimes b_{ij}}$, we have
%
\begin{eqnarray}
&& \bigl\langle\bigl(D^i(\sqrt{t} F_j + \sqrt{1-t}
S_j\eta_j)\bigr)^{\otimes
b_{ij}} , v \bigr
\rangle_{\EuFrak H^{\otimes b_{ij} }}
\nonumber\\
&&\label{e2}\qquad = \sum_{u=0}^{b_{ij}} \pmatrix{b_{ij}
\cr
u} t^{u/2} (1-t)^{(b_{ij}- u)/2} \eta^{(b_{ij}- u)}\\
&&\hspace*{6pt} \qquad \qquad {}\times \bigl\langle
\bigl(D^iF_j\bigr)^{\otimes u} \otimes
\bigl(D^iS_j\bigr)^{\otimes(b_{ij} -u)} , v \bigr
\rangle_{\EuFrak H^{\otimes
b_{ij}}}.\nonumber
\end{eqnarray}
Substituting (\ref{e2})  into (\ref{e1}), and taking into account the
symmetry of $u_k$, yields
\begin{eqnarray*}
&& E \bigl[ \bigl\langle D^{q_k}\varphi_{x_k}(\mathbf{X} ;
\sqrt{t}
F+\sqrt{1-t}S\cdot\eta),u_k \bigr
\rangle_{\mathfrak{H}^{\otimes
q_k}} \bigr]
\\
&&\qquad = \sum_{\beta^{(q_k)} \in\mathscr{B}(q_k)} C\bigl(\alpha^{(q_k)}\bigr)
t^{|b'|/2} (1-t)^{|b''|/2} \prod_{i=1}^{q_k}
\prod_{j=1}^d \pmatrix{b'_{ij}+b''_{ij}
\cr
b'_{ij}}
\\
&& \quad \qquad {}\times E\Biggl[\partial^{\alpha(\beta
^{(q_k)})}\varphi_{x_k} (\mathbf{X} ;
\sqrt{t}%
F+\sqrt{1-t}S\cdot\eta)\prod_{j=1}^d
\eta_j^{|b''_{\bullet j }|}
\\
&& \qquad \quad{}\times \Biggl\langle u_k, h_1^{\otimes a_1}
\otimes\cdots
\otimes h_m^{\otimes a_m} \bigotimes
_{i=1}^{q_k} \bigotimes_{j=1}^{d}
\bigl\{ \bigl(D^iF_j\bigr)^{\otimes
b'_{ij}}\otimes
\bigl(D^iS_j\bigr)^{\otimes b''_{ij}} \bigr\} \Biggr\rangle
_{\EuFrak H
^{\otimes q_k}}\Biggr]\!.
\end{eqnarray*}

Notice that if $\beta^{(q_k)} $ does not belong to $\mathscr
{B}_0(q_k)$, then $b'_{q_kl}\ge1$ for some index $l=1,\ldots, d$.
Taking into account the relations (\ref{Dio}) this implies that
$b'_{q_kl}= 1$, $b'_{q_kj}= 0$ for all $j\neq l$, $k_{q_k}=1$ and all
the other entries of $\beta^{(q_k)} $ must be equal to zero. In this
way, the above sum can be decomposed as follows:
\begin{eqnarray*}
&& \sum_{\beta^{(q_k)} \in\mathscr{B}_0(q_k)} C\bigl(\alpha^{(q_k)}\bigr)
t^{|b'|/2} (1-t)^{|b''|/2} \prod_{i=1}^{q_k}
\prod_{j=1}^d \pmatrix{b'_{ij}+b''_{ij}
\cr
b'_{ij}}
\\
&&\quad{}\times E\Biggl[\partial^{\alpha(\beta
^{(q_k)})}\varphi_{x_k} (\mathbf{X} ;
\sqrt{t}%
F+\sqrt{1-t}S\cdot\eta)\prod_{j=1}^d
\eta_j^{|b''_{\bullet j }|}
\\
&&\quad {}\times \Biggl\langle u_k, h_1^{\otimes a_1}
\otimes\cdots
\otimes h_m^{\otimes a_m} \bigotimes
_{i=1}^{q_k} \bigotimes_{j=1}^{d}
\bigl\{ \bigl(D^iF_j\bigr)^{\otimes
b'_{ij}}\otimes
\bigl(D^iS_j\bigr)^{\otimes b''_{ij}} \bigr\} \Biggr\rangle
_{\EuFrak H
^{\otimes q_k}} \Biggr]
\\
&&\quad {}+ \sum_{l=1} ^d \sqrt{t} E \biggl[
\frac{\partial^2}{\partial x_k\,
\partial x_l}\varphi(\mathbf{X} ; \sqrt{t}%
F+\sqrt{1-t}S\cdot\eta)
\bigl\langle D^{q_k}F_l , u_k\bigr
\rangle_{{\EuFrak
H^{\otimes
q_k}}} \biggr]
\\
&&\qquad:= D(k,t)+F(k,t).
\end{eqnarray*}
Since
\[
\Biggl\llvert \frac{%
1}{2\sqrt{t}}\sum_{k=1}^d
F(k,t) - \frac{%
1}{2}\sum_{k=1}^d E
\biggl[\frac{\partial^2}{\partial x_k^2} \varphi (\mathbf{X} ; \sqrt{t}F+\sqrt{1-t}S\cdot
\eta)S_k^{2}%
 \biggr]\Biggr\rrvert
\leqslant(\ref{w1}),
\]
the theorem is proved once we show that
\[
\sum_{k=1}^d \int_0^1
\frac{1}{2\sqrt{t}} \bigl|D(k,t) \bigr| \,dt
\]
is less than the sum in (\ref{w2}). Using the independence of $\eta$
and $X$, conditioning with respect to $X$ and applying Lemma~\ref
{lmipp} yields
\begin{eqnarray*}
&& E\Biggl[\partial^{\alpha(\beta^{(q_k)})}\varphi_{x_k} (\mathbf{X} ;
\sqrt{t}%
F+\sqrt{1-t}S\cdot\eta)\prod_{j=1}^d
\eta_j^{|b''_{\bullet j }|}
\\
&& \quad{}\times \Biggl\langle u_k, h_1^{\otimes
a_1}
\otimes\cdots\otimes h_m^{\otimes a_m} \bigotimes
_{i=1}^{q_k} \bigotimes_{j=1}^{d}
\bigl\{ \bigl(D^iF_j\bigr)^{\otimes b'_{ij}}\otimes
\bigl(D^iS_j\bigr)^{\otimes b''_{ij}} \bigr\} \Biggr
\rangle_{\EuFrak
H^{\otimes q_k}} \Biggr]
\\
&&\qquad = \sum_{l_1=0}^{\lfloor|b''_{\bullet1}|/2\rfloor} \cdots \sum
_{l_d=0}^{\lfloor|b''_{\bullet d}|/2\rfloor} \prod_{s=1}^d
\frac{|
b''_{\bullet s}| !}{2^{l_s}(| b''_{\bullet s}| - 2l_s)!l_s!}
\\
&& \quad\qquad{}\times E\Biggl[ \Biggl\langle u_k, h_1^{\otimes a_1}
\otimes\cdots \otimes h_m^{\otimes a_m} \bigotimes
_{i=1}^{q_k} \bigotimes_{j=1}^{d}
\bigl\{ \bigl(D^iF_j\bigr)^{\otimes b'_{ij}}\otimes
\bigl(D^iS_j\bigr)^{\otimes
b''_{ij}} \bigr\} \Biggr
\rangle_{\EuFrak H^{\otimes q_k}}
\\
&&\hspace*{35pt}\hspace*{26pt} \quad\qquad{}\times \prod_{s=1}^d
S^{|b''_{\bullet s}| -2l_s} \partial_\star^{(\beta
^{(q_k)};l_1,\ldots,l_d)}\varphi_{x_k}
(\mathbf{X} ; \sqrt{t}F+\sqrt{1-t}S\cdot\eta)^{} \Biggr].
\end{eqnarray*}

Then, estimating the term
$\partial_\star^{(\beta^{(q_k)};l_1,\ldots,l_d)}\varphi_{x_k} (\mathbf{X} ;
\sqrt{t}F+\sqrt{1-t}S\cdot\eta)$ by $\|\partial_\star^{(\beta
^{(q_k)};l_1,\ldots,l_d)}\varphi_{x_k}\|_\infty$, which does not depend on
$t$, and using the equation
\[
\int_0^1 \frac{1}{\sqrt{t} } t^{|b'|/2}
(1-t) ^{ |b''|/2} \,dt= B\bigl( |b'|/2 +1/2,
|b''|/2+1\bigr),
\]
we obtain the desired estimate.
\end{pf}

\section{Application to weighted quadratic variations}\label{sex2}

In this section, we apply the previous results to the case of weighted
quadratic variations of fractional Brownian motion.
Let us introduce first some notation.

Given a measurable function $f\dvtx \mathbb{R}\rightarrow\mathbb{R}$, an
integer $N\ge0$
and a real number $p\ge1$ we define the seminorm
%
\begin{equation}
\label{eq3a}
\|f\|_{N,p} = \sum_{i=0}^N
\sup_{0\leqslant t\leqslant1} \bigl\| f^{(i)} \bigr\|_{
L^p(\mathbb{R}, \gamma_t)},
\end{equation}
where $\gamma_t$ is\vadjust{\goodbreak} the normal distribution $N(0,t)$.

We say that a function $f\dvtx \mathbb{R}\rightarrow\mathbb{R}$ has
\textit{moderate
growth} if there exist positive constants $A$, $B$ and $\alpha<2$ such that
for all $x \in\mathbb{R}$, $|f(x)| \leqslant A \exp ( B |x|
^\alpha )$.
Notice that the seminorm (\ref{eq3a}) is finite if $f$ and all its
derivatives up to the order $N$ have moderate growth.

Consider a fractional Brownian motion $B=\{B_t\dvtx  t\in[0,1]\}$ with
Hurst parameter $H\in(0,1)$.
That is, $B$ is a zero mean Gaussian process with covariance
$ E(B_tB_s) =\frac{1}2  ( t^{2H} + s^{2H} - |t-s| ^{2H}  )$.
The process $B$ can be extended to an isonormal Gaussian process
indexed by the Hilbert space $\EuFrak{H}$, which is the closure of the
set of simple functions on $[0,1]$ with respect to the inner product
$\langle\mathbf{1}_{[0,t]} , \mathbf{1}_{[0,s]} \rangle_{\EuFrak
{H}} =
E(B_tB_s)$. We refer the reader to the basic references \cite{nour,nualartbook} for a detailed account on this process. We denote by
%
\begin{equation}
\label{rho}
\rho_H(k)= \tfrac{1}2 \bigl(
|k+1|^{2H} + |k-1|^{2H} -2 |k|^{2H} \bigr),\qquad k\in
\mathbb{Z},
\end{equation}
the covariance function of the stationary sequence $\{ B(k+1)-B(k)\dvtx
k\ge0\}$.

We consider the uniform partition of the interval $[0,1]$, and for any
$n\ge1$ and $k=0,\ldots, n-1$ we denote $\Delta
B_{k/n}=B_{(k+1)/n}-B_{k/n}$, $\delta_{k/n}=\mathbf{1}_{[k/n,(k+1)/n]}$
and~$\varepsilon_{k,n}=\mathbf{1}_{[0, k/n]}$. We will also make use of the notation
$\beta_{j,k}=\langle\delta_{j/n} ,\delta_{k/n} \rangle_{\EuFrak H }$
and $\alpha_{j,t}=\langle\delta_{j/n} ,\mathbf{1}_{[0,t]} \rangle
_{\EuFrak H }$, for any $t\in[0,1]$ and $j,k=0,\ldots, n-1$.

Given a function $f\dvtx \mathbb{R}\rightarrow\mathbb{R}$, we define
\[
u_n= n^{2H-{1}/2} \sum_{k=0}^{n-1}
f(B_{k/n}) \delta_{k/n}^{\otimes2}.
\]
We are interested in the asymptotic behavior of the weighted quadratic
functionals
%
\begin{eqnarray}
F_n &=&  n^{2H-{1}/2} \sum
_{k=0}^{n-1} f(B_{k/n}) \bigl[ (\Delta
B_{k/n})^2 - n^{-2H} \bigr]
\nonumber
\\[-8pt]
\label{eq10}\\[-8pt]
\nonumber
&=& n^{2H-{1}/2}
\sum_{k=0}^{n-1} f(B_{k/n})I_2
\bigl(\delta _{k/n}^{\otimes2}\bigr).
\end{eqnarray}
It is known (see, e.g., \cite{nou,nounua,nounuatud}) that for
$\frac{1}4 <H<\frac{3}4$, $F_n$ converges in law to a mixture of Gaussian
distributions. When the Hurst parameter $H$ is not in this range, a
different phenomenon occurs, as it was observed by Nourdin in \cite{nou}. More precisely, for $H<\frac{1}4$, $n^{2H-{1}/2} F_n$
converges in $L^2(\Omega)$ to $\frac{1}4 \int_0^1 f''(B_s)\,ds$, whereas
for $H>\frac{3}4$,
$n^{{3}/2 -2H} F_n$ converges in $L^2(\Omega)$ to $\int_0^1 f(B_s)
\,dZ_s$, where $Z$ is the Rosenblatt process
(see \cite{nounuatud,nou}).
In the critical case $H=\frac{1}4$, there is convergence in law to a
linear combination of the limits in the cases $H<\frac{1}4$ and $\frac{1}4<H< \frac{3}4$, and in the critical case $H=\frac{3}4$ there is
convergence in law with an additional logarithmic factor (see \cite{nounuatud,nou}).

In view of these results, we will focus on the case $\frac{1}4<H<\frac{3}4$, although our result could easily be extended to the limit case
$H=\frac{3}4$. Outside the interval $ [ \frac{1}4, \frac{3}4 ]$
the convergence is in $L^2(\Omega)$ and our methodology does not seem
to be well suited to study the rate of convergence.
Applying the general approach developed in previous sections, we are
able to show the following rate of convergence in the asymptotic
behavior of $F_n$, in the case $H\in (\frac{1}4, \frac{3}4 )$.
This represents a quantitative version of the convergence in law proved
in \cite{nounuatud}.

\begin{prop}\label{prop61}
Assume that the Hurst index $H$ of $B$ belongs to $(\frac{1}4,\frac{3}4)$.
Consider a function $f\dvtx \mathbb{R}\to\mathbb{R}$ of class $C^4$ such
that $f$ and its
first $4$ derivatives have moderate growth.
Suppose in addition that\break $E  [  ( \int_0^1 f^2(B_s)\,ds  )
^{-\alpha}  ] <\infty$ for some $\alpha>1$.
Consider the sequence of random variables $F_n$ defined by (\ref
{eq10}). Set $S=\sqrt{\sigma_H \int_0^1 f^2(B_s)\,ds}$, with $\sigma^2_H=
\sum_{k=-\infty}^{\infty} \rho_H(k)^2$, where $\rho_H$ is defined in
(\ref{rho}). Then, for any function $\varphi\dvtx \mathbb{R}\rightarrow
\mathbb{R}$
of class $C^{5}$ with $\Vert\varphi^{(k)}\Vert_{\infty}<\infty$
for any $%
k=0,\ldots,5$ we have
%
\begin{equation}
\label{eq66}
\quad\bigl|E\bigl[\varphi(F_n)\bigr]-E\bigl[\varphi( S\eta)
\bigr]\bigr| \leqslant C_{f,H} \max_{1\leqslant i\leqslant5} \bigl\|
\varphi^{(i)} \bigr\|_\infty n^{ -( |2H-{1}/2| \wedge|2H-{3}/2| ) },
\end{equation}
where $\eta$ is a standard normal variable independent of $B$. The
constant $C_{f,H}$ has the form $ C_{f,H}= C_H \max ( 1, \|f\|
_{4,4}^4,  (1+ |E[S^{-2\alpha}] |^{{1}/\alpha} \|f\|
_{1,5\beta}
^5  ) )$, where $C_H$ depends on $H$ and $\frac{1}\alpha+
\frac{1}\beta=1$.
\end{prop}

\begin{pf}
Along the proof $C$ will
denote a generic constant that might depend on $H$.

Notice first that the random variable $F_n$ does not coincide with
$\delta^2(u_n)$, except in the case $H=\frac{1}2$. For this reason, we define
$G_n= \delta^2(u_n)$,
and show the following estimate for the difference $F_n-G_n$:
%
\begin{equation}\label{eq32}
E\bigl[|F_n-G_n|\bigr] \leqslant C\|f\|_{3,2}
n^{-( |2H-{1}/2| \wedge
|2H-{3}/2| )}.
\end{equation}
To show (\ref{eq32}), we first apply Lemma~\ref{lemme} and we obtain
\[
F_n-G_n = n^{2H-{1}/2} \sum
_{k=0}^{n-1} 2\delta \bigl( f'(B_{k/n})
\delta_{k/n} \bigr) \alpha_{k,k/n} + n^{2H-{1}/2} \sum
_{k=0}^{n-1} f''(B_{k/n})
\alpha_{k,k/n}^2.
\]
Using the equality $\delta ( f'(B_{k/n}) \delta_{k/n}  ) = f'(B_{k/n})
I_1(\delta
_{k/n} )-f''(B_{k/n}) \alpha_{k,k/n}$,
yields
\begin{eqnarray*}
F_n-G_n &=&  2 n^{2H-{1}/2} \sum
_{k=0}^{n-1} f'(B_{k/n})
I_1(\delta _{k/n} )\alpha_{k,k/n} -
n^{2H-{1}/2} \sum_{k=0}^{n-1}
f''(B_{k/n}) \alpha_{k,k/n}^2
\\
&\hspace*{-3pt}:=&  2M_n-R_n.
\end{eqnarray*}
Point (a) of Lemma~\ref{lem1} implies $|\alpha_{k,k/n}| \le
n^{-(2H)\wedge1}$ and we can write
%
\begin{equation}\label{eq30}
E\bigl[ |R_n|\bigr] \leqslant\|f\|_{2,1} n^{{1}/2+2H -(4H\wedge2)}.
\end{equation}
On the other hand,
\[
E\bigl[M_n^2\bigr] = n^{4H-1} \sum
_{j,k=0}^{n-1} E\bigl[f'(B_{j/n})
f'(B_{k/n}) I_1(\delta_{j/n}
)I_1(\delta_{k/n} )\bigr] \alpha_{j,j/n}
\alpha_{k,k/n},
\]
and using the relation
\[
I_1(\delta_{j/n} )I_1(\delta_{k/n}
)= I_2( \delta_{j/n} \,\widetilde{\otimes}\,
\delta_{k/n} ) + \langle\delta_{j/n}, \delta_{k/n}
\rangle_{\EuFrak H}
\]
the duality relationship (\ref{ipp}) yields
\begin{eqnarray*}
E\bigl[M_n^2\bigr] &\leqslant & \|f\|_{3,2}^2
n^{4H-1} \sum_{j,k=0}^{n-1} \bigl[| \beta
_{j,k} |+|\alpha_{j,j/n} \alpha_{k,k/n}|+|
\alpha_{j,k/n} \alpha _{k,j/n}|\bigr]\\
&&\hspace*{60pt}\qquad {}\times | \alpha_{j,j/n}
\alpha_{k,k/n}|.
\end{eqnarray*}
Finally, applying points (a) and (c) of Lemma~\ref{lem1}, we obtain,
%
\begin{equation}
\label{eq31}
E\bigl[M_n^2\bigr] \leqslant C\|f
\|_{3,2}^2 n^{4H-1} \bigl(n^{(1-2H)\vee0} +
n^{2 -(4H\wedge2)} \bigr) n^{-(4H\wedge2)}.
\end{equation}
If $H <\frac{1}2$, we obtain a rate of the form $n^{1-4H}$ and if
$H\geqslant
\frac{1}2$ we obtain the bound
$n^{4H-3}$. Then the estimates (\ref{eq30}) and (\ref{eq31})
imply (\ref{eq32}).

Taking into account the estimate (\ref{eq32}), the estimate (\ref
{eq66}) will follow from (\ref{eq11}), provided we show the following
inequalities for some constant $C$ depending on $H$ and for any $\beta>1$:
%
\begin{eqnarray}
\label{1}
E \bigl( \bigl|\bigl\langle u_n, D^2G_n\bigr\rangle_{{\EuFrak H}^{\otimes2}} - S^2 \bigr| \bigr)&\leqslant& C\|f
\|^2_{4,2} n^{-( |2H-{1}/2| \wedge
|2H-{3}/2| ) },
\\
\label{2}
E \bigl( \bigl|\bigl\langle u_n, DG_n^{\otimes2}\bigr
\rangle_{{\EuFrak
H}^{\otimes
2}} \bigr| \bigr)&\leqslant& C \|f\|^3_{3,3}
n^{-( |2H-{1}/2|
\wedge
|2H-{3}/2| ) },
\\
\label{3}
\bigl\llVert \bigl\langle u_n, D\bigl(S^2
\bigr)^{\otimes2}\bigr\rangle_{{\EuFrak H}^{\otimes
2}}\bigr\rrVert _{ L^\beta(\Omega)} &
\leqslant& C\| f\|^5_{1, 5\beta} n^{-(
|2H-{1}/2| \wedge|2H-{3}/2|)},
\\
\label{4}
E \bigl( \bigl|\bigl\langle u_n, D^2\bigl(S^2
\bigr)\bigr\rangle_{{\EuFrak H}^{\otimes2}} \bigr| \bigr)&\leqslant& C\|f\|^3_{2,3}
n^{ -( |2H-{1}/2| \wedge
|2H-{3}/2| ) },
\\
\label{5}
E \bigl( \bigl|\bigl\langle u_n, DG_n\otimes D
\bigl(S^2\bigr)\bigr\rangle_{{\EuFrak
H}^{\otimes2}} \bigr| \bigr)&\leqslant& C \|f
\|^4_{3,4} n^{-(
|2H-{1}/2| \wedge|2H-{3}/2|)}.
\end{eqnarray}
The derivatives $S^2$ are given by the following expressions:
\begin{eqnarray*}
D\bigl(S^2\bigr)&=& 2\sigma_H\int_0^1
\bigl(ff'\bigr) (B_s)\mathbf{1}_{[0,s]}\,ds,
\\[-2pt]
D^2\bigl(S^2\bigr)&=& 2\sigma_H \int
_0^1 \bigl(f'^2+ff''
\bigr) (B_s)\mathbf{1}_{[0,s]^2}\,ds.
\end{eqnarray*}
On the other hand, applying formula (\ref{t5}) we obtain the following
expressions for the derivatives of $G_n$
\begin{eqnarray*}
DG_n&=& \delta(u_n)+ \delta^2
(Du_n),
\\[-2pt]
D^2G_n&=& u_n +2 \delta(Du_n)
+ \delta^2\bigl( D^2 u_n\bigr).
\end{eqnarray*}
We are now ready to prove (\ref{1})--(\ref{5}). The proof will be based
on the estimates obtained in Lemma~\ref{lem2} of the \hyperref[sec7]{Appendix}.
\end{pf}

\begin{pf*}{Proof of (\ref{1})}
We have
\begin{eqnarray*}
&& \bigl\llvert \bigl\langle u_n, D^2G_n\bigr
\rangle_{\EuFrak H^{\otimes2}} - S^2\bigr\rrvert \\[-2pt]
&& \qquad \leqslant  \bigl\llvert
\|u_n\|^2_{\EuFrak H^{\otimes2}} - S^2 \bigr\rrvert
+ 2 \bigl\llvert \bigl\langle u_n, \delta(Du_n) \bigr
\rangle_{\EuFrak H^{\otimes2}} \bigr\rrvert + \bigl\llvert \bigl\langle u_n,
\delta^2\bigl( D^2 u_n\bigr)\bigr
\rangle_{\EuFrak H^{\otimes
2}} \bigr\rrvert
\\[-2pt]
&& \qquad =: |A_n|+ 2 |B_n| + |C_n|.
\end{eqnarray*}
To estimate $E[|A_n|]$, we write
\begin{eqnarray*}
\|u_n\|^2_{\EuFrak H^{\otimes2}}& =&n^{4H-1}\sum
_{j,k=0}^{n-1} f(B_{j/n})f(B_{k/n})
\beta_{j,k}^{2} \\[-2pt]
&=& \frac{1}{ n}\sum
_{j,k=0}^{n-1} f(B_{j/n})f(B_{k/n})
\rho_H(k-j)^{2}
\\[-2pt]
&=&\frac{1}{ n}\sum_{p=-n+1}^{n-1}\sum
_{j=0\vee-p}^{ (
n-1 ) \wedge ( n-1-p ) } f(B_{j/n})f(B_{ (
j+p ) /n})
\rho_H(p) ^{2}.
\end{eqnarray*}
If we replace $f(B_{ (j+p ) /n})$ by $ f(B_{j/n})$ we make an
error in expectation of $(p/n)^H$, so this produces a total error of
$n^{-H}$. On the other hand, the sequence
$\sum_{|p|>n} \rho_H(p) ^{2}$
converges to zero at the rate $n^{4H-3}$.
As a consequence,
\begin{eqnarray*}
E\bigl[|A_n| \bigr]  &\leqslant &  C \bigl(\|f\|^2_{1,2}
n^{-H}+ \|f\|^2_{0,2} n^{4H-3}\bigr)\\
&&{}+
\sigma _H^2 E \Biggl[ \Biggl\llvert \frac{1}n
\sum_{k=0}^{n-1} f^2(B_{k/n})
- \int_0^1 f^2(B_s)\,ds
\Biggr\rrvert \Biggr].
\end{eqnarray*}
It remains to estimate
\[
\frac{1}n \sum_{k=0}^{n-1}
f^2(B_{k/n}) - \int_0^1
f^2(B_s)\,ds = \sum_{k=0}^{n-1}
\int_{k/n}^{(k+1)/n} \bigl[f^2(B_{k/n})-f^2(B_s)
\bigr]\,ds.
\]
Using that $E[|f^2(B_{k/n})-f^2(B_s)|] \leqslant C \|f\| ^2_{1,2} n^{-H}$
for $s\in[k/n,(k+1)/n]$, we obtain:
%
\begin{equation}
\label{x11}
E\bigl[|A_n| \bigr] \leqslant C \bigl(\|f\|^2_{1,2}
n^{-H}+ \|f\|^2_{0,2} n^{4H-3}\bigr).
\end{equation}
For the term $B_n$ we can write, using (\ref{m2}) and Meyer's inequalities:
%
\begin{eqnarray}
E\bigl[|B_n|\bigr]  &\leqslant &  n^{2H-{1}/2} \sum
_{k=0}^{n-1} E\bigl[\bigl| f(B_{k/n}) \delta
\bigl( D _{k/n} (u_n \otimes_1
\delta_{k/n})\bigr)\bigr|\bigr]
\nonumber
\\[-8pt]
\label{x22}\\[-8pt]
\nonumber
& \leqslant &  C\|f\|^2_{2,2}
n^{2H-(3H\wedge{3}/2)}.
\end{eqnarray}
The term $C_n$ is handled in the same way, by using Meyer's
inequalities and point~(d) of Lemma~\ref{lem1}:
%
\begin{eqnarray}
E\bigl[|C_n|\bigr] &\leqslant& n^{2H-{1}/2} \sum
_{k=0}^{n-1} E\bigl[\bigl| f(B_{k/n}) \delta
^2\bigl( D^2 _{k/n,k/n} u_n \bigr)\bigr|
\bigr]
\nonumber
\\[-8pt]
\label{x3}
\\[-8pt]
\nonumber
 &\leqslant& C n^{2H-{1}/2} \|f\|^2_{4,2}
\sum_{j,k=0}^{n-1} \beta_{j,k}^2
\le C n^{{1}/2-2H} \|f\|^2_{4,2}.
\end{eqnarray}
Then (\ref{1}) follows from (\ref{x11}), (\ref{x22}) and (\ref{x3}).
\end{pf*}

\begin{pf*}{Proof of (\ref{2})}
We have
\begin{eqnarray*}
\bigl\langle u_n, DG_n^{\otimes2}\bigr
\rangle_{{\EuFrak H}^{\otimes2} } &=& \bigl\langle u_n, \delta(u_n)
\otimes\delta(u_n) \bigr\rangle_{{\EuFrak
H}^{\otimes2} }+ 2 \bigl\langle
u_n, \delta(u_n) \otimes\delta^2
(Du_n) \bigr\rangle_{{\EuFrak
H}^{\otimes2} }
\\
&& {}+ \bigl\langle u_n, \delta^2 (Du_n)
\otimes\delta^2 (Du_n) \bigr\rangle _{{\EuFrak
H}^{\otimes2} }
\\
&=:& A_n+ 2B_n + C_n.
\end{eqnarray*}
For the term $A_n$ we have, applying H\"older's and Meyer's
inequalities and the estimate (\ref{m1}),
\begin{eqnarray*}
E\bigl[|A_n| \bigr] &\leqslant &  n^{2H-{1}/2} \sum
_{k=0}^{n-1} E\bigl[ \bigl|f(B_{k/n}) \bigl(
\delta (u_n\otimes_ 1\delta_{k/n})
\bigr)^2 \bigr|\bigr] \\
&\leqslant &  C \|f\|^3_{1,3}
n^{2H-{1}/2-
(2H\wedge1)}.
\end{eqnarray*}
Similarly, using H\"older's and Meyer's inequalities and the estimates
(\ref{m1}) and~(\ref{m3}) yields
\begin{eqnarray*}
E\bigl[|B_n|\bigr] &\leqslant &  n^{2H-{1}/2} \sum
_{k=0}^{n-1} E\bigl[ \bigl|f(B_{k/n})
\delta(u_n\otimes_
1\delta_{k/n})
\delta^2(D_{k/n}u_n)\bigr|\bigr]\\
& \leqslant &  C \|f
\|^3_{3,3} n^{2H-
(3H\wedge {3}/2)}.
\end{eqnarray*}
Finally, using again H\"older's and Meyer's inequalities and the
estimate (\ref{m1}) yields
\begin{eqnarray*}
E\bigl[|C_n|\bigr] & \leqslant &  n^{2H-{1}/2} \sum
_{k=0}^{n-1} E\bigl[ \bigl|f(B_{k/n}) \bigl(
\delta ^2(D_{k/n}u_n)\bigr)^2\bigr|
\bigr] \\
& \leqslant &  C \|f\|^3_{3,3} n^{2H-{1}/2- (2H\wedge1)}.
\end{eqnarray*}
\upqed\end{pf*}

\begin{pf*}{Proof of (\ref{3})}
We have
\begin{eqnarray*}
&& \bigl\langle u_n, D\bigl(S^2\bigr)^{\otimes2}\bigr
\rangle_{{\EuFrak H}^{\otimes2} }\\
&& \qquad  =16n^{2H-{1}/2}\sum_{k=0}^{n-1}
f(B_{k/n}) \int_0^1 \!\int
_0^1 \bigl(ff'\bigr)
(B_s) \bigl(ff'\bigr) (B_t)
\alpha_{k,t} \alpha_{k,s}\, ds\,dt.
\end{eqnarray*}
Then, we can write, using points (a) and (b) of Lemma~\ref{lem1},
\begin{eqnarray*}
E \bigl[ \bigl\llvert \bigl\langle u_n, D\bigl(S^2
\bigr)^{\otimes2}\bigr\rangle_{{\EuFrak
H}^{\otimes2} } \bigr\rrvert \bigr]  &\leqslant &  C
\| f\| ^5_{1,5\beta} n^{2H-{1}/2} \sup
_{s,t \in
[0,1]}\sum_{k=0}^{n-1} |
\alpha_{k,t} \alpha_{k,s} | \\
&\leqslant &  C \|f\|
^5_{1,5\beta}n^{2H-{1}/2- (2H\wedge1)},
\end{eqnarray*}
for any $\beta\geqslant1$.
\end{pf*}

\begin{pf*}{Proof of (\ref{4})}
We have
\[
\bigl\langle u_n, D^2\bigl(S^2\bigr)\bigr
\rangle_{{\EuFrak H}^{\otimes2} } = 4n^{2H-{1}/2}\sum_{k=0}^{n-1}
f(B_{k/n})\int_0^1
\bigl(f'^2+ff''\bigr)
(B_s) \alpha_{k,t}^2\,ds.
\]
As a consequence, applying points (a) and (b) of Lemma~\ref{lem1} yields
\begin{eqnarray*}
E \bigl[ \bigl\llvert \bigl\langle u_n, D^2
\bigl(S^2\bigr)\bigr\rangle_{{\EuFrak H}^{\otimes2} } \bigr\rrvert \bigr]
& \leqslant &  C \|f\|^3_{2,3} n^{2H-{1}/2} \sup
_{s \in[0,1]}\sum_{k=0}^{n-1}
\alpha_{k,s}^2 \\
&\leqslant &  C\|f\|^3_{2,3}
n^{2H-{1}/2- (2H\wedge1)}.
\end{eqnarray*}
\upqed\end{pf*}

\begin{pf*}{Proof of (\ref{5})}
We have
\begin{eqnarray*}
\bigl\langle u_n, DG_n\otimes D\bigl(S^2
\bigr)\bigr\rangle_{{\EuFrak H}^{\otimes2} } &=&  \bigl\langle u_n,
\delta(u_n)\otimes D\bigl(S^2\bigr)\bigr
\rangle_{{\EuFrak H}^{\otimes2}
} + \bigl\langle u_n, \delta^2(Du_n)
\otimes D\bigl(S^2\bigr)\bigr\rangle_{{\EuFrak
H}^{\otimes
2} }
\\
&\,\,=:& A_n+B_n.
\end{eqnarray*}
For the term $A_n$ we can write, applying H\"older's and Meyer's
inequalities and the estimate (\ref{m1}),
\begin{eqnarray*}
E\bigl[|A_n|\bigr] & \leqslant & n^{2H-{1}/2} \sum
_{k=0}^{n-1} E\bigl[\bigl| f(B_{k/n}) \delta
(u_n\otimes_1\delta_{k/n} ) D_{k/n}
\bigl(S^2\bigr)\bigr|\bigr]\\
& \leqslant &  C \|f\|^4_{1,4}
n^{2H-{1}/2- (2H\wedge1)}.
\end{eqnarray*}
For the term $A_n$ we can write, applying H\"older's and Meyer's
inequalities and the estimate (\ref{m1}),
\begin{eqnarray*}
E\bigl[|A_n|\bigr]  &\leqslant &  n^{2H-{1}/2} \sum
_{k=0}^{n-1} E\bigl[\bigl| f(B_{k/n}) \delta
^2(D_{k/n}u_n ) D_{k/n}
\bigl(S^2\bigr)\bigr|\bigr] \\
&\leqslant &  C \|f\|^4_{3,4}
n^{2H+{1}/2- (4H\wedge2)}.
\end{eqnarray*}
This completes the proof of Proposition~\ref{prop61}.
\end{pf*}

%
\begin{rem}
Note that the exponent in the rate $\delta= -(|2H-\frac{1}2|
\wedge|2H-\frac{3}2| )$ is minimum when $H=\frac{1}2 $ with $\delta
=-\frac{1}2$. On the other hand, it becomes worst when $H$ goes away
from $\frac{1}2$ either from below or from above, and it converges to
zero as $H$ tends to $\frac{1}4$ or $\frac{3}4$. This is natural in view
of the limit results for the weighted quadratic variations obtained in
\cite{nou,nounuatud}.
This phenomenon has not been observed in other asymptotic problems,
such as the rate of convergence for Euler-type numerical approximations
of stochastic differential equations, where the rate $-(2H-\frac{1}2)$
improves when $H$ increases from $\frac{1}2$ up to $\frac{3}4$ (see \cite{huliunua}).
\end{rem}

\begin{rem}
In the case $H=\frac{1}2$, the process $B$ is a Brownian motion,
and it has independent increments. As  consequence $\beta_{j,k} =0$
for $j\neq k$. Moreover, $F_n=G_n$. Therefore, the estimate (\ref{eq66}) can be replaced by
\[
\bigl|E\bigl[\varphi(F_n)\bigr]-E\bigl[\varphi( S\eta)\bigr]\bigr| \leqslant
C_{f } \max_{2\leqslant i\leqslant5} \bigl\|\varphi^{(i)}
\bigr\|_\infty n^{ - {1}/2 },
\]
where $S^2 =2\int_0^1 f(B_s)^2 \,ds$.
\end{rem}

\begin{rem}
The extension to weighted power variations of any order or to
Euler numerical schemes for stochastic differential equations driven by
a fractional Brownian motion seems more involved. In the case of Euler
numerical schemes, the results that could be obtained applying the
methodology developed in this paper would lead to a precise analysis of
the rate of convergence of the error to a particular distribution,
which is usually a mixture of Gaussian laws. That is, we would be able
to establish how close is the error to a limit distribution in terms of
a distance between probabilities defined by means of regular\vadjust{\goodbreak} functions.
\end{rem}

\begin{appendix}
\section*{Appendix}\label{sec7}

In this section, we will show two technical lemmas that play a
fundamental role in the analysis of the asymptotic quadratic variation
of the fractional Brownian motion. The notation in both lemmas is taken
from Section~\ref{sex2}.

\setcounter{prop}{0}
\begin{lemme} \label{lem1}
Let $0<H<1$ and $n\geqslant1$. We have, for some
constant $C_H$:
\begin{longlist}[(d)]
\item[(a)] $\llvert  \alpha_{k,t} \rrvert  \leqslant n^{-(2H\wedge1)}$
for any $t\in
{}[0,1]$ and $k=0,\ldots, n-1$.
\item[(b)] \vspace*{1pt}$\sup_{t\in{}[0,1]}\sum_{k=0}^{n-1}\llvert  \alpha
_{k,t}\rrvert  \leqslant C_H$.
\item[(c)] $\sum_{k,j=0}^{n-1}\llvert  \beta_{j,k}\rrvert  \leqslant C_H
n^{(1-2H)\vee0}$.
\item[(d)] If $H<\frac{3}4$, then $\sum_{k,j=0}^{n-1} \beta_{j,k}
^{2}\leqslant C_Hn^{1-4H}$.
\item[(e)] $\sum_{k,j=0}^{n-1}\llvert  \beta_{k,l} \beta_{j,l}
\rrvert
\leqslant C_{H}n^{-(4H \wedge2)}$ for any $l=0,\ldots, n-1$.
\item[(f)] If $H<\frac{3}4$, then $\sum_{k,j=0}^{n-1}\llvert  \beta_{k,l},
\beta_{j,l} \beta_{j,k}\rrvert  \leqslant C_H n^{-4H -(2H\wedge1)}$
for any
$l=0,\ldots, n-1$.
\end{longlist}
\end{lemme}

\begin{pf}
Parts (a), (c) and (d) are contained in Lemmas 5 and 6 of \cite{nounuatud}. Part~(b) has been proved in Lemma~5.1 of \cite{nounua} in
the case $H<\frac{1}2$ and the proof actually works for any $H\in
(0,1)$. Part (e) follows easily from
\[
\sum_{k,j=0}^{n-1}\llvert
\beta_{k,l} \beta_{j,l} \rrvert = \frac{1}4
n^{-4H}\sum_{k,j=0}^{n-1}\bigl|
\rho_H(k-l) \rho_H(j-l)\bigr|,
\]
and the fact that the series $\sum_{p\in\mathbb{Z}} |\rho_H(p)|$ is
convergent if $0<H\leqslant\frac{1}2$ and it diverges at the rate $n^{2H-1}$
if $H >\frac{1}2$. Finally, to prove (f) we write, using Young's inequality,
\begin{eqnarray*}
\sum_{k,j=0}^{n-1}\llvert
\beta_{k,l} \beta_{j,l} \beta_{j,k}\rrvert &=&
\frac{1}8 n^{-6H}\sum_{k,j=0}^{n-1}\bigl|
\rho_H(k-l) \rho_H(j-l)\rho _H(j-k)\bigr|
\\
&\leqslant& \frac{1}8 n^{-6H} \biggl( \sum
_{p\in\mathbb{Z}} \rho _H(p)^2 \biggr) \Biggl(
\sum^n_{p=-n} \bigl|\rho_H(p)\bigr|
\Biggr)
\\
&\leqslant& C_H n^{-4H -(2H\wedge1)},
\end{eqnarray*}
where we have exploited the fact that $ \sum_{p\in\mathbb{Z}} \rho
_H(p)^2$ is convergent (because $H<\frac{3}4$), together with the
asymptotic behavior of the mapping $n\mapsto\break \sum^n_{p=-n} |\rho_H(p)|$.
\end{pf}

The next lemma provides some technical estimates.

\begin{lemme} \label{lem2}
For any integer $M\geqslant0$ and any real
number $p>1$, there exists a constant $C$ depending on $M, p$ and the
Hurst parameter $H$ such that:
%
\setcounter{equation}{0}
\begin{eqnarray}\label{m1}
\| u_n \otimes_1 \delta_{k/n}
\|_{M,p} & \leqslant& C \|f\|_{M,p} n^{
-{1}/2- ( H\wedge{1}/2)},
\\
\label{m2}
\bigl\| D_{k/n} (u_n \otimes_1
\delta_{k/n} )\bigr\|_{M,p} &\leqslant& C \|f\| _{M+1,p}
n^{ -{1}/2 - ( 3H\wedge {3}/2) },
\\
\label{m3}
\| D_{k/n} u_n \|_{M,p} &
\leqslant& C \|f\|_{M+1,p} n^{-(2H\wedge1)},
\end{eqnarray}
where $D_{k/n} F$ means $\langle DF, \delta_{k/n}\rangle_{\EuFrak
H}$, for a
given random variable $F$.
\end{lemme}

\begin{pf}
In order to show the first estimate, we can write, for any integer
$0\leqslant m \leqslant M$,
\[
D^m( u_n \otimes_1 \delta_{k/n})=n^{2H-{1}/2}
\sum_{j=0}^{n-1} f^{(m)}(B_{j/n})
\beta_{k,j} \delta_{j/n}\, \widetilde{\otimes}\, \varepsilon
_{j/n}^{\otimes m}.
\]
Then, using points (a), (e) and (f) of Lemma~\ref{lem1} we obtain
\begin{eqnarray*}
&& \bigl(E \bigl[ \bigl\| D^m( u_n \otimes_1
\delta_{k/n}) \bigr\|^p_{\EuFrak H
^{\otimes(m+1)}} \bigr]
\bigr)^{{1}/p} \\
&& \qquad \leqslant Cn^{2H-{1}/2} \|f\|_{m,p}
\\
&& \qquad \quad {}\times \Biggl( \sum_{j,j'=0}^{n-1} \bigl|
\beta_{k,j} \beta_{k,j'} \bigl\langle\delta_{j/n}
\,\widetilde{\otimes}\, \varepsilon_{j/n}^{\otimes m},
\delta_{j'/n} \,\widetilde{\otimes}\, \varepsilon_{j'/n}^{\otimes m}
\bigr\rangle _{\EuFrak H^{\otimes(m+1)}}\bigr| \Biggr)^{{1}/2}
\\
&& \qquad\leqslant C \|f\|_{m,p} n^{2H-{1}/2} \Biggl(\sum
_{j,j'=0}^{n-1}|\beta_{k,j}
\beta_{k,j'}| \bigl( |\beta_{j,j'} | +|\alpha_{j,j'/n}
\alpha_{j',j/n}| \bigr) \Biggr)^{{1}/2}
\\
&& \qquad \leqslant C \|f\|_{m,p} n^{2H-{1}/2 } \bigl( n^{-2H -(H\wedge
{1}/2 )}
+ n^{-(4H\wedge2)} \bigr)\\
&& \qquad \leqslant C \|f\|_{m,p} n^{ -{1}/2- ( H\wedge{1}/2)},
\end{eqnarray*}
which shows (\ref{m1}).

To show the second estimate, we can write, for any integer $0\leqslant m
\leqslant M$,
\[
D^m D_{k/n} (u_n \otimes_1
\delta_{k/n} )=n^{2H-{1}/2} \sum_{j=0}^{n-1}
f^{(m+1)}(B_{j/n}) \beta_{k,j} \alpha_{k,j/n}
\delta_{j/n} \,\widetilde{\otimes}\, \varepsilon_{j/n}^{\otimes m}.
\]
Then, using points (a), (e) and (f) of Lemma~\ref{lem1} we obtain
\begin{eqnarray*}
&& \bigl(E \bigl[ \bigl\| D^m D_{k/n} (u_n
\otimes_1 \delta_{k/n} )\bigr\| ^p_{\EuFrak H
^{\otimes(m+1)}}
\bigr] \bigr)^{{1}/p}\\
&& \qquad  \leqslant Cn^{2H-{1}/2} \|f\|_{m+1,p}
\\
&&\qquad \quad {}\times \Biggl( \sum_{j,j'=0}^{n-1} \bigl|
\beta_{k,j} \alpha_{k,j/n} \beta_{k,j'}
\alpha_{k,j'/n} \bigl\langle\delta_{j/n} \,\widetilde {\otimes}\,
\varepsilon_{j/n}^{\otimes m}, \delta_{j'/n} \,\widetilde{
\otimes} \,\varepsilon _{j'/n}^{\otimes m} \bigr\rangle_{\EuFrak H^{\otimes(m+1)}}\bigr|
\Biggr)^{{1}/2}
\\
&& \qquad\leqslant Cn^{ 2H-{1}/2- ( 2H\wedge1) } \|f\|_{m+1,p}\\
&&\qquad \quad  {}\times \Biggl(\sum
_{j,j'=0}^{n-1}|\beta_{k,j}
\beta_{k,j'}|\bigl( |\beta _{j,j'}|+|\alpha _{j, j'/n}
\alpha_{j',j/n} |\bigr) \Biggr)^{{1}/2}
\\
&& \qquad\leqslant C \|f\|_{m+1,p}n^{ 2H-{1}/2- ( 2H\wedge1) } \bigl(n^{-2H-(H\wedge{1}/2)} +
n^{-(4H\wedge2)} \bigr)
\\
&& \qquad\leqslant C \|f\|_{m+1,p}n^{ -{1}/2- ( 3H\wedge{3}/2)},
\end{eqnarray*}
and (\ref{m2}) follows.

Finally,
for the estimate (\ref{m3}) we can write
\[
D^m D_{k/n} u_n =n^{2H-{1}/2}\sum
_{j=0}^{n-1} f^{(m+1)}(B_{j/n})
\alpha_{k,j/n} \delta_{j/n}^{\otimes2} \,\widetilde{\otimes}\,
\varepsilon _{j/n}^{\otimes m},
\]
which implies, using points (a), (c) and (d) of Lemma~\ref{lem1},
\begin{eqnarray*}
&& \bigl(E \bigl[ \bigl\| D^m D_{k/n} u_n
\bigr\|^p_{\EuFrak H^{\otimes(m+2)}} \bigr] \bigr)^{{1}/p} \\
&& \qquad \leqslant
Cn^{2H-{1}/2} \|f\|_{m+1,p}
\\
&&\quad  \qquad {}\times \Biggl( \sum_{j,j'=0}^{n-1} \bigl|
\alpha_{k,j/n} \alpha _{k,j'/n} \bigl\langle\delta_{j/n}
^{\otimes2}\, \widetilde{\otimes} \,\varepsilon _{j/n}^{\otimes m},
\delta_{j'/n} ^{\otimes2} \,\widetilde{\otimes}\, \varepsilon_{j'/n}^{\otimes m}
\bigr\rangle_{\EuFrak H^{\otimes(m+2)}}\bigr| \Biggr)^{{1}/2}
\\
&& \qquad \leqslant Cn^{2H-{1}/2-(2H\wedge1) } \|f\|_{m+1,p}\\
&& \qquad \quad {}\times\Biggl( \sum
_{j,j'=0}^{n-1} \bigl(\beta_{j,j'}^2
+|\beta_{j,j'} \alpha_{j,j'/n}\alpha _{j',j/n}| +
\alpha_{j,j'/n}^2 \alpha_{j',j/n}^2 \bigr)
\Biggr)^{{1}/2}
\\
&& \qquad \leqslant C \|f\|_{m+1,p} n^{2H-{1}/2-(2H\wedge1) } \\
&& \qquad  \quad{}\times\bigl( n^{{1}/2-2H}
+ n^{ [({1}/2-H) \vee0] -(2H\wedge1)} +n^{ -(4H\wedge2) } \bigr).
\end{eqnarray*}
This shows (\ref{m3}) and the proof of the lemma is complete.
\end{pf}
\end{appendix}

\section*{Acknowledgments}
We are grateful to two anonymous referees for a
thorough reading
and a number of helpful suggestions.






\printaddresses

\begin{thebibliography}{42}


\bibitem{BNCP}
\begin{barticle}[mr]
\bauthor{\bsnm{Barndorff-Nielsen},~\bfnm{Ole~E.}\binits{O.~E.}},
\bauthor{\bsnm{Corcuera},~\bfnm{Jos{\'e}~Manuel}\binits{J.~M.}} \AND
\bauthor{\bsnm{Podolskij},~\bfnm{Mark}\binits{M.}}
(\byear{2009}).
\btitle{Power variation for {G}aussian processes with stationary increments}.
\bjournal{Stochastic Process. Appl.}
\bvolume{119}
\bpages{1845--1865}.
\bid{doi={10.1016/j.spa.2008.09.004}, issn={0304-4149}, mr={2519347}}
\end{barticle}
%

\bptok{imsref}%
\endbibitem

%

\bibitem{boupec}
\begin{barticle}[mr]
\bauthor{\bsnm{Bourguin},~\bfnm{Solesne}\binits{S.}} \AND
\bauthor{\bsnm{Peccati},~\bfnm{Giovanni}\binits{G.}}
(\byear{2014}).
\btitle{Portmanteau inequalities on the {P}oisson space: Mixed regimes and multidimensional clustering}.
\bjournal{Electron. J. Probab.}
\bvolume{19}
\bpages{no. 66, 1--42}.
\bid{doi={10.1214/EJP.v19-2879}, issn={1083-6489}, mr={3248195}}
\bptnote{check volume, check pages, check year}%
\end{barticle}
%
%

\bptok{imsref}%
\endbibitem

%

\bibitem{chatterjee}
\begin{barticle}[mr]
\bauthor{\bsnm{Chatterjee},~\bfnm{Sourav}\binits{S.}}
(\byear{2010}).
\btitle{Spin glasses and {S}tein's method}.
\bjournal{Probab. Theory Related Fields}
\bvolume{148}
\bpages{567--600}.
\bid{doi={10.1007/s00440-009-0240-8}, issn={0178-8051}, mr={2678899}}
\end{barticle}
%

\bptok{imsref}%
\endbibitem

%

\bibitem{cnw}
\begin{barticle}[mr]
\bauthor{\bsnm{Corcuera},~\bfnm{Jos{\'e}~Manuel}\binits{J.~M.}},
\bauthor{\bsnm{Nualart},~\bfnm{David}\binits{D.}} \AND
\bauthor{\bsnm{Woerner},~\bfnm{Jeannette~H.~C.}\binits{J.~H.~C.}}
(\byear{2006}).
\btitle{Power variation of some integral fractional processes}.
\bjournal{Bernoulli}
\bvolume{12}
\bpages{713--735}.
\bid{doi={10.3150/bj/1155735933}, issn={1350-7265}, mr={2248234}}
\end{barticle}
%

\bptok{imsref}%
\endbibitem

\bibitem{D}
\begin{bmisc}[auto:parserefs-M02]
\bauthor{\bsnm{D{\"{o}}bler},~\bfnm{C.}\binits{C.}}
(\byear{2013}).
\bhowpublished{Stein's method of exchangeable pairs for absolutely continuous, univariate distributions with applications to the Polya urn model.
Preprint.}
\end{bmisc}
%

\bptok{imsref}%
\endbibitem

\bibitem{EdVie}
\begin{bincollection}[auto:parserefs-M02]
\bauthor{\bsnm{Eden},~\bfnm{R.}\binits{R.}} \AND
\bauthor{\bsnm{Viens},~\bfnm{F.}\binits{F.}}
(\byear{2013}).
\btitle{General upper and lower tail estimates using Malliavin calculus and Stein's equations}.
In \bbooktitle{Seminar on Stochastic Analysis, Random Fields and Applications VII}
(\beditor{\bfnm{R.~C.}\binits{R.~C.}~\bsnm{Dalang}},
\beditor{\bfnm{M.}\binits{M.}~\bsnm{Dozzi}} \AND
\beditor{\bfnm{F.}\binits{F.}~\bsnm{Russo}}, eds.).
\bseries{Progress in Probability}
\bvolume{67}.
\bpublisher{Birkh\"{a}user},
\blocation{Basel}.
\end{bincollection}
%

\bptok{imsref}%
\endbibitem

\bibitem{EdViq}
\begin{barticle}[mr]
\bauthor{\bsnm{Eden},~\bfnm{R.}\binits{R.}} \AND
\bauthor{\bsnm{V\'{\i}quez},~\bfnm{J.}\binits{J.}}
(\byear{2012}).
\btitle{Nourdin--Peccati analysis on Wiener and Wiener--Poisson space for general distributions}.
\bjournal{Stochastic Process. Appl.}
\bjournal{125}
\bpages{182--216}.
\bid{mr={3274696}}
\end{barticle}
%

\bptok{imsref}%
\endbibitem

\bibitem{HarNua2}
\begin{barticle}[mr]
\bauthor{\bsnm{Harnett},~\bfnm{Daniel}\binits{D.}} \AND
\bauthor{\bsnm{Nualart},~\bfnm{David}\binits{D.}}
(\byear{2012}).
\btitle{Weak convergence of the {S}tratonovich integral with respect to a class of {G}aussian processes}.
\bjournal{Stochastic Process. Appl.}
\bvolume{122}
\bpages{3460--3505}.
\bid{doi={10.1016/j.spa.2012.06.008}, issn={0304-4149}, mr={2956113}}
\end{barticle}
%

\bptok{imsref}%
\endbibitem

\bibitem{HarNua1}
\begin{barticle}[mr]
\bauthor{\bsnm{Harnett},~\bfnm{Daniel}\binits{D.}} \AND
\bauthor{\bsnm{Nualart},~\bfnm{David}\binits{D.}}
(\byear{2013}).
\btitle{Central limit theorem for a {S}tratonovich integral with {M}alliavin calculus}.
\bjournal{Ann. Probab.}
\bvolume{41}
\bpages{2820--2879}.
\bid{doi={10.1214/12-AOP769}, issn={0091-1798}, mr={3112933}}
\end{barticle}
%

\bptok{imsref}%
\endbibitem

\bibitem{huliunua}
\begin{bmisc}[auto:parserefs-M02]
\bauthor{\bsnm{Hu},~\bfnm{Y.}\binits{Y.}},
\bauthor{\bsnm{Liu},~\bfnm{Y.}\binits{Y.}} \AND
\bauthor{\bsnm{Nualart},~\bfnm{D.}\binits{D.}}
(\byear{2013}).
\bhowpublished{Modified Euler approximation scheme for stochastic differential equations driven by fractional Brownian motions.
Preprint.}
\end{bmisc}
%

\bptok{imsref}%
\endbibitem

\bibitem{JacSh}
\begin{bbook}[mr]
\bauthor{\bsnm{Jacod},~\bfnm{Jean}\binits{J.}} \AND
\bauthor{\bsnm{Shiryaev},~\bfnm{Albert~N.}\binits{A.~N.}}
(\byear{1987}).
\btitle{Limit Theorems for Stochastic Processes}.
\bseries{Grundlehren der Mathematischen Wissenschaften [Fundamental Principles of Mathematical Sciences]}
\bvolume{288}.
\bpublisher{Springer},
\blocation{Berlin}.
\bid{doi={10.1007/978-3-662-02514-7}, mr={0959133}}
\end{bbook}
%

\bptok{imsref}%
\endbibitem

\bibitem{KT1}
\begin{barticle}[mr]
\bauthor{\bsnm{Kusuoka},~\bfnm{Seiichiro}\binits{S.}} \AND
\bauthor{\bsnm{Tudor},~\bfnm{Ciprian~A.}\binits{C.~A.}}
(\byear{2012}).
\btitle{Stein's method for invariant measures of diffusions via {M}alliavin calculus}.
\bjournal{Stochastic Process. Appl.}
\bvolume{122}
\bpages{1627--1651}.
\bid{doi={10.1016/j.spa.2012.02.005}, issn={0304-4149}, mr={2914766}}
\end{barticle}
%

\bptok{imsref}%
\endbibitem

\bibitem{KT2}
\begin{bmisc}[auto:parserefs-M02]
\bauthor{\bsnm{Kusuoka},~\bfnm{S.}\binits{S.}} \AND
\bauthor{\bsnm{Tudor},~\bfnm{C.~A.}\binits{C.~A.}}
(\byear{2013}).
\bhowpublished{Extension of the Fourth Moment Theorem to invariant measures of diffusions.
Preprint.}
\end{bmisc}
%

\bptok{imsref}%
\endbibitem

\bibitem{Mis}
\begin{barticle}[mr]
\bauthor{\bsnm{Mishkov},~\bfnm{Rumen~L.}\binits{R.~L.}}
(\byear{2000}).
\btitle{Generalization of the formula of {F}aa di {B}runo for a composite function with a vector argument}.
\bjournal{Int. J. Math. Math. Sci.}
\bvolume{24}
\bpages{481--491}.
\bid{doi={10.1155/S0161171200002970}, issn={0161-1712}, mr={1781515}}
\end{barticle}
%

\bptok{imsref}%
\endbibitem

\bibitem{nou}
\begin{barticle}[mr]
\bauthor{\bsnm{Nourdin},~\bfnm{Ivan}\binits{I.}}
(\byear{2008}).
\btitle{Asymptotic behavior of weighted quadratic and cubic variations of fractional {B}rownian motion}.
\bjournal{Ann. Probab.}
\bvolume{36}
\bpages{2159--2175}.
\bid{doi={10.1214/07-AOP385}, issn={0091-1798}, mr={2478679}}
\end{barticle}
%

\bptok{imsref}%
\endbibitem

\bibitem{nour}
\begin{bbook}[mr]
\bauthor{\bsnm{Nourdin},~\bfnm{Ivan}\binits{I.}}
(\byear{2012}).
\btitle{Selected Aspects of Fractional {B}rownian Motion}.
\bseries{Bocconi \& Springer Series}
\bvolume{4}.
\bpublisher{Springer},
\blocation{Milan}.
\bid{doi={10.1007/978-88-470-2823-4}, mr={3076266}}
\end{bbook}
%

\bptok{imsref}%
\endbibitem

\bibitem{nounua}
\begin{barticle}[mr]
\bauthor{\bsnm{Nourdin},~\bfnm{Ivan}\binits{I.}} \AND
\bauthor{\bsnm{Nualart},~\bfnm{David}\binits{D.}}
(\byear{2010}).
\btitle{Central limit theorems for multiple {S}korokhod integrals}.
\bjournal{J. Theoret. Probab.}
\bvolume{23}
\bpages{39--64}.
\bid{doi={10.1007/s10959-009-0258-y}, issn={0894-9840}, mr={2591903}}
\end{barticle}
%

\bptok{imsref}%
\endbibitem

\bibitem{nounuatud}
\begin{barticle}[mr]
\bauthor{\bsnm{Nourdin},~\bfnm{Ivan}\binits{I.}},
\bauthor{\bsnm{Nualart},~\bfnm{David}\binits{D.}} \AND
\bauthor{\bsnm{Tudor},~\bfnm{Ciprian~A.}\binits{C.~A.}}
(\byear{2010}).
\btitle{Central and non-central limit theorems for weighted power variations of fractional {B}rownian motion}.
\bjournal{Ann. Inst. Henri Poincar\'{e} Probab. Stat.}
\bvolume{46}
\bpages{1055--1079}.
\bid{doi={10.1214/09-AIHP342}, issn={0246-0203}, mr={2744886}}
\end{barticle}
%

\bptok{imsref}%
\endbibitem

\bibitem{noupecwei}
\begin{barticle}[mr]
\bauthor{\bsnm{Nourdin},~\bfnm{Ivan}\binits{I.}} \AND
\bauthor{\bsnm{Peccati},~\bfnm{Giovanni}\binits{G.}}
(\byear{2008}).
\btitle{Weighted power variations of iterated {B}rownian motion}.
\bjournal{Electron. J. Probab.}
\bvolume{13}
\bpages{1229--1256}.
\bid{doi={10.1214/EJP.v13-534}, issn={1083-6489}, mr={2430706}}
\end{barticle}
%

\bptok{imsref}%
\endbibitem

%

%

\bibitem{np-book}
\begin{bbook}[mr]
\bauthor{\bsnm{Nourdin},~\bfnm{Ivan}\binits{I.}} \AND
\bauthor{\bsnm{Peccati},~\bfnm{Giovanni}\binits{G.}}
(\byear{2012}).
\btitle{Normal Approximations with {M}alliavin Calculus: From Stein's Method to Universality}.
\bseries{Cambridge Tracts in Mathematics}
\bvolume{192}.
\bpublisher{Cambridge Univ. Press},
\blocation{Cambridge}.
\bid{doi={10.1017/CBO9781139084659}, mr={2962301}}
\end{bbook}
%

\bptok{imsref}%
\endbibitem

%

%

\bibitem{nou-poly}
\begin{barticle}[mr]
\bauthor{\bsnm{Nourdin},~\bfnm{Ivan}\binits{I.}} \AND
\bauthor{\bsnm{Poly},~\bfnm{Guillaume}\binits{G.}}
(\byear{2013}).
\btitle{Convergence in total variation on {W}iener chaos}.
\bjournal{Stochastic Process. Appl.}
\bvolume{123}
\bpages{651--674}.
\bid{doi={10.1016/j.spa.2012.10.004}, issn={0304-4149}, mr={3003367}}
\end{barticle}
%

\bptok{imsref}%
\endbibitem

\bibitem{nourev}
\begin{barticle}[mr]
\bauthor{\bsnm{Nourdin},~\bfnm{Ivan}\binits{I.}} \AND
\bauthor{\bsnm{R{\'e}veillac},~\bfnm{Anthony}\binits{A.}}
(\byear{2009}).
\btitle{Asymptotic behavior of weighted quadratic variations of fractional {B}rownian motion: The critical case {$H=1/4$}}.
\bjournal{Ann. Probab.}
\bvolume{37}
\bpages{2200--2230}.
\bid{doi={10.1214/09-AOP473}, issn={0091-1798}, mr={2573556}}
\end{barticle}
%

\bptok{imsref}%
\endbibitem

\bibitem{nourevswa}
\begin{barticle}[mr]
\bauthor{\bsnm{Nourdin},~\bfnm{Ivan}\binits{I.}},
\bauthor{\bsnm{R{\'e}veillac},~\bfnm{Anthony}\binits{A.}} \AND
\bauthor{\bsnm{Swanson},~\bfnm{Jason}\binits{J.}}
(\byear{2010}).
\btitle{The weak {S}tratonovich integral with respect to fractional {B}rownian motion with {H}urst parameter {$1/6$}}.
\bjournal{Electron. J. Probab.}
\bvolume{15}
\bpages{2117--2162}.
\bid{doi={10.1214/EJP.v15-843}, issn={1083-6489}, mr={2745728}}
\end{barticle}
%

\bptok{imsref}%
\endbibitem

\bibitem{nualartbook}
\begin{bbook}[mr]
\bauthor{\bsnm{Nualart},~\bfnm{David}\binits{D.}}
(\byear{2006}).
\btitle{The {M}alliavin Calculus and Related Topics},
\bedition{2nd} ed.
\bseries{Probability and Its Applications (New York)}.
\bpublisher{Springer},
\blocation{Berlin}.
\bid{mr={2200233}}
\end{bbook}
%

\bptok{imsref}%
\endbibitem

\bibitem{nuaort}
\begin{barticle}[mr]
\bauthor{\bsnm{Nualart},~\bfnm{D.}\binits{D.}} \AND
\bauthor{\bsnm{Ortiz-Latorre},~\bfnm{S.}\binits{S.}}
(\byear{2008}).
\btitle{Central limit theorems for multiple stochastic integrals and {M}alliavin calculus}.
\bjournal{Stochastic Process. Appl.}
\bvolume{118}
\bpages{614--628}.
\bid{doi={10.1016/j.spa.2007.05.004}, issn={0304-4149}, mr={2394845}}
\end{barticle}
%

\bptok{imsref}%
\endbibitem

\bibitem{nunugio}
\begin{barticle}[mr]
\bauthor{\bsnm{Nualart},~\bfnm{David}\binits{D.}} \AND
\bauthor{\bsnm{Peccati},~\bfnm{Giovanni}\binits{G.}}
(\byear{2005}).
\btitle{Central limit theorems for sequences of multiple stochastic integrals}.
\bjournal{Ann. Probab.}
\bvolume{33}
\bpages{177--193}.
\bid{doi={10.1214/009117904000000621}, issn={0091-1798}, mr={2118863}}
\end{barticle}
%

\bptok{imsref}%
\endbibitem

\bibitem{petamwii}
\begin{barticle}[mr]
\bauthor{\bsnm{Peccati},~\bfnm{Giovanni}\binits{G.}} \AND
\bauthor{\bsnm{Taqqu},~\bfnm{Murad~S.}\binits{M.~S.}}
(\byear{2008}).
\btitle{Stable convergence of multiple {W}iener--{I}t\^o integrals}.
\bjournal{J. Theoret. Probab.}
\bvolume{21}
\bpages{527--570}.
\bid{doi={10.1007/s10959-008-0154-x}, issn={0894-9840}, mr={2425357}}
\end{barticle}
%

\bptok{imsref}%
\endbibitem

\bibitem{pectud}
\begin{bincollection}[mr]
\bauthor{\bsnm{Peccati},~\bfnm{Giovanni}\binits{G.}} \AND
\bauthor{\bsnm{Tudor},~\bfnm{Ciprian~A.}\binits{C.~A.}}
(\byear{2005}).
\btitle{Gaussian limits for vector-valued multiple stochastic integrals}.
In \bbooktitle{S\'eminaire de {P}robabilit\'es {XXXVIII}}.
\bseries{Lecture Notes in Math.}
\bvolume{1857}
\bpages{247--262}.
\bpublisher{Springer},
\blocation{Berlin}.
\bid{mr={2126978}}
\bptnote{check year}%
\end{bincollection}
%

\bptok{imsref}%
\endbibitem

\bibitem{py1}
\begin{bincollection}[mr]
\bauthor{\bsnm{Peccati},~\bfnm{Giovanni}\binits{G.}} \AND
\bauthor{\bsnm{Yor},~\bfnm{Marc}\binits{M.}}
(\byear{2004}).
\btitle{Four limit theorems for quadratic functionals of {B}rownian motion and {B}rownian bridge}.
In \bbooktitle{Asymptotic Methods in Stochastics}.
\bseries{Fields Inst. Commun.}
\bvolume{44}
\bpages{75--87}.
\bpublisher{Amer. Math. Soc.},
\blocation{Providence, RI}.
\bid{mr={2106849}}
\end{bincollection}
%

\bptok{imsref}%
\endbibitem

\bibitem{py2}
\begin{bincollection}[mr]
\bauthor{\bsnm{Peccati},~\bfnm{Giovanni}\binits{G.}} \AND
\bauthor{\bsnm{Yor},~\bfnm{Marc}\binits{M.}}
(\byear{2004}).
\btitle{Hardy's inequality in {$L\sp 2([0,1])$} and principal values of {B}rownian local times}.
In \bbooktitle{Asymptotic Methods in Stochastics}.
\bseries{Fields Inst. Commun.}
\bvolume{44}
\bpages{49--74}.
\bpublisher{Amer. Math. Soc.},
\blocation{Providence, RI}.
\bid{mr={2106848}}
\end{bincollection}
%

\bptok{imsref}%
\endbibitem

\bibitem{podvetter}
\begin{barticle}[mr]
\bauthor{\bsnm{Podolskij},~\bfnm{Mark}\binits{M.}} \AND
\bauthor{\bsnm{Vetter},~\bfnm{Mathias}\binits{M.}}
(\byear{2010}).
\btitle{Understanding limit theorems for semimartingales: A short survey}.
\bjournal{Stat. Neerl.}
\bvolume{64}
\bpages{329--351}.
\bid{doi={10.1111/j.1467-9574.2010.00460.x}, issn={0039-0402}, mr={2683464}}
\end{barticle}
%

\bptok{imsref}%
\endbibitem

\bibitem{R}
\begin{bincollection}[mr]
\bauthor{\bsnm{Reinert},~\bfnm{Gesine}\binits{G.}}
(\byear{2005}).
\btitle{Three general approaches to {S}tein's method}.
In \bbooktitle{An Introduction to {S}tein's Method}.
\bseries{Lect. Notes Ser. Inst. Math. Sci. Natl. Univ. Singap.}
\bvolume{4}
\bpages{183--221}.
\bpublisher{Singapore Univ. Press},
\blocation{Singapore}.
\bid{doi={10.1142/9789812567680_0004}, mr={2235451}}
\end{bincollection}
%

\bptok{imsref}%
\endbibitem

\bibitem{reny}
\begin{barticle}[mr]
\bauthor{\bsnm{R{\'e}nyi},~\bfnm{Alfr{\'e}d}\binits{A.}}
(\byear{1963}).
\btitle{On stable sequences of events}.
\bjournal{Sankhya, Ser. A}
\bvolume{25}
\bpages{293--302}.
\bid{issn={0581-572X}, mr={0170385}}
\end{barticle}
%

\bptok{imsref}%
\endbibitem

\bibitem{rev2009}
\begin{barticle}[mr]
\bauthor{\bsnm{R{\'e}veillac},~\bfnm{Anthony}\binits{A.}}
(\byear{2009}).
\btitle{Convergence of finite-dimensional laws of the weighted quadratic variations process for some fractional {B}rownian sheets}.
\bjournal{Stoch. Anal. Appl.}
\bvolume{27}
\bpages{51--73}.
\bid{doi={10.1080/07362990802564491}, issn={0736-2994}, mr={2473140}}
\end{barticle}
%

\bptok{imsref}%
\endbibitem

\bibitem{RY}
\begin{bbook}[mr]
\bauthor{\bsnm{Revuz},~\bfnm{Daniel}\binits{D.}} \AND
\bauthor{\bsnm{Yor},~\bfnm{Marc}\binits{M.}}
(\byear{1999}).
\btitle{Continuous Martingales and {B}rownian Motion},
\bedition{3rd} ed.
\bseries{Grundlehren der Mathematischen Wissenschaften [Fundamental Principles of Mathematical Sciences]}
\bvolume{293}.
\bpublisher{Springer},
\blocation{Berlin}.
\bid{doi={10.1007/978-3-662-06400-9}, mr={1725357}}
\end{bbook}
%

\bptok{imsref}%
\endbibitem
\end{thebibliography}
\end{document}